\documentstyle [12pt,twoside, epsf]{article}
 \newcommand{\ZZ}{\mbox{$Z\!\!\! Z\!$}} 

 \let\\=\cr
 
 \def\Chi{\hbox{\raise0.5ex\hbox{$\chi$}}}

 \newtheorem{th}{Theorem}
 \newtheorem{lem}{Lemma}
 
 \newtheorem{prop}{Proposition}

 \newtheorem{defn}{Definition}
 \newtheorem{conj}{Conjecture}

 \newtheorem{rem}{Remark}

\def\picill#1by#2(#3)
{\vbox to #2
{\hrule width #1 height 0pt depth 0pt
\vfill\epsffile{#3}}}

\begin{document}
 \pagestyle{myheadings}

 \markboth{{\sc Kauffman \& Lambropoulou}}{{\sc  Hard Unknots and Collapsing Tangles}}

 \title{\bf Hard Unknots and Collapsing Tangles}

\author{Louis H. Kauffman
 and
Sofia Lambropoulou 
}

\date{}

\maketitle


\section{Introduction}
This paper gives infinitely many examples of unknot diagrams that are {\it hard,} in the sense that the diagrams
need to be made more complicated by Reidemeister moves before they can be simplified. In order to construct these
diagrams, we prove theorems characterizing when the numerator of the sum of two rational tangles is an unknot.
The paper uses these results in studying processive DNA recombination, finding minimal size unknot diagrams,
generalizing to collapses to knots as well as to unknots, and in finding unknots with arbirarily high complexity in terms of 
the {\it recalcitrance} that will be defined below.
\bigbreak

We have made every effort to make this paper self-contained and available to a reader who is just becoming interested
in the theory of knots and links. We believe that a good way to appreciate the problem of knotting is to look at the structure of
unknots. The fact that our problems are elucidated by the theory of rational tangles (see below) provides a good opportunity
for introducing the basics of this theory and its intimate connection with matrices and continued fractions. Because the paper is
both a research paper and an expository paper, we have designed the rest of this introduction and the next section as a guide,
and the following section as a  review of the needed tangle theory.
\bigbreak 

The fundamental problem in knot theory is to determine whether a closed loop embedded in three dimensional
space is knotted. We sometimes put this as the question: Is this knot knotted? This might seem to be a 
contradiction in terms, but that comes from referring to any tangled loop as a ``knot", when only some loops
are so irretrievably tangled that there is no way to simplify them. View Figure 6. Can you tell whether this knot is 
knotted or not? It requires quite some intuition for topology to just look at a knot and know if it is 
really knotted.
\bigbreak

In order to analyze the knottedness of a knot, a mathematical representation is required.
In this paper we shall use the method of knot and link diagrams, and the equivalence relation generated by the 
Reidemeister moves (See Figure 1 for an illustration of these moves). Knot diagrams are graphs with extra structure that 
encode the embedding type of the knot. Each diagram is a pictorial representation of the knot, and so appeals to the 
intuition of the viewer. The Reidemeister moves are a set of simple combinatorial moves, proved in the 
1920's to capture the notion of topological equivalence of knots and links in three dimensional space. Single applications of 
these moves can leave the diagram with the same number of crossings (places where a weaving of two segments occurs), or 
increase or decrease the number of crossings. Some unknottings can be accomplished without increasing the number of crossings
in the diagram. We call such unknot diagrams {\it easy} since the fact that they are unknotted can be determined by a 
finite search for simplifying moves. However, there are culprit diagrams that require moves that {\it increase} the 
number of crossings before the diagram can be simplified to an unknotted circle with no crossings. It is the structure of such 
culprits that is the  subject matter of this paper.
\bigbreak

See Figure 2  for a diagram that we shall refer throughout this paper as the ``Culprit." This culprit is not the only 
culprit, but it is the exemplar that we shall use, and it is the example that started this investigation. The first 
author likes to use the Culprit as an example in introductory talks about knot theory. One draws the Culprit on the board
and asks whether it is knotted or not. This gives rise to a discussion of easy and hard unknots, and how the existence of 
hard unknots makes us need a theory of knots in order to prove knottedness when it occurs. After using this example, we
began to ask how to produce other examples that were hard and to wonder if our familiar culprit might be the  smallest such
example (size being the number of crossings, in this case 10).
\bigbreak 

We discovered that there are infinitely many 
examples of hard unknot diagrams, obtained by using the theory, due to
John H. Conway, of rational tangles and their closures 
We discuss the clues to how this comes about in Section 1 of this paper, and we give other examples of culprits and near
culprits in this section. In order to use the theory of rational tangles, one must become familiar with the notion of {\it tangle}
and the notion of the {\it fraction of a tangle.} In Section 1, we introduce the tangle analysis and assume that the reader
knows about tangle fractions. The reader can skip to Section 2 to learn this material, and then finish reading 
Section 1, if this is needed. Section 2 contains the basics about the theory of rational tangles, and their closures,
the rational knots and links. 
\bigbreak

Section 1 contains our definition of the {\it recalcitrance} $R(D)$ of an unknot diagram $D$. The recalcitrance of 
$D$ is the ratio of the least number of crossings in a possibly larger diagram that is needed in an unknotting sequence for $D,$
to the number of crossings in $D.$ The recalcitrance is a measure of how much complexification is needed in order to 
simplify $D.$
\bigbreak

The clue in Section 1 is the fact that the Culprit can be divided into two rational tangles 
whose fractions add up to a fraction whose numerator has absolute value equal to 1. It turns out that {\it whenever the sum of the fractions of
two rational tangles has numerator equal to plus or minus one, then the closure of the sum of the two tangles will be an unknot}.
This result is Theorem 5 in Section 3.  In Theorem 8 of Section 4 we
take a further step and characterize fractions
$\frac{P}{Q}$ and $\frac{R}{S}$ such that $\frac{P}{Q} - \frac{R}{S} = \frac{\pm 1}{QS}$ in terms of their associated 
continued fractions. It turns out that this last equation is satisfied if and only if one of the two continued fractions is a 
{\it convergent} of the other. This means that one continued fraction is a one-term truncate of the other. For example
$$3/2 = 1 + \frac{1}{2}$$ is a convergent of $$10/7 = 1 + \frac{1}{2 + \frac{1}{3}}.$$ 
Section 2
sets up the matrix representations for continued fractions that underpin the proof of the Theorems. This completely solves the
question of when two fractions give rise to an unknot via the (numerator) closure of the sum of their associated tangles.
\bigbreak

In Section 5 we use these results to construct many examples of hard unknots. The first example, $K=N([1,4]-[1,3]),$ of this section is given
in Figure 24 and its mirror image $H$ in Figure 31. This culprit $K$ is a hard unknot diagram with only 9 crossings.
We then show how our original 
culprit (of 10 crossings) arises from a ``tucking construction" applied to an unknot that is an easy diagram
without the tuck (Figure 26). This section then discusses other applications of the tucking construct. In Section 6 we prove
that the $9$ crossing examples of Figure 31 and some relatives obtained by flyping and taking
mirror images (see Figure 13 for the defintion of ``flype") are the smallest hard unknot diagrams that can be made by taking the closure of
the sum of two alternating rational tangles. In Section 7 we show the historically first hard unknot, due to Goeritz in 1934. The Goeritz
diagram has 11 crossings, but there is a surprise: In unknotting the Goeritz diagram, we find the culprit $H$ appearing (Figure
33).
\bigbreak

In Section 8, we show how our unknots are related to the study of processive recombination of DNA.
In the {\it tangle model} for DNA recombination, pioneered by DeWitt Sumners and Claus Ernst, the initial substrate
of the DNA is represented as the closure of the sum of two rational tangles. It is usual to assume that the initial 
DNA substrate is unknotted. We have characterized such unknot configurations in this paper, and so are in a position to 
apply our results to the model. We show that processive recombination stabilizes, in the sense that 
the form of the resulting knotted or linked DNA is obtained by just adding twists {\it in a single site} on the closure of
a certain tangle. This result helps to understand the form of the recombination process.
\bigbreak

In Section 9 we consider more general collapses of diagrams to knots and links rather than just to unknots.
Thus we generalize our results to hard knots as well as hard unknots.
In Section 10 we revisit the recalcitrance of unknot diagrams, and show, by using tangle methods,
that there is no upper bound to the 
recalcitrance. 
\bigbreak

Section 11 is a vignette of other mathematical relationships associated with rational numbers and continued fractions.
In particular, we show how the Farey series gives rise to a natural way to find pairs of fractions $\frac{a}{b}$ and $\frac{c}{d}$ such that 
$|ad-bc| = 1,$ and how the continued fraction expansions of all rational numbers and real numbers arise in a context
that already involves geometry and topology.  All these relationships add
dimension to our discussion of unknots related to such fraction pairs. By looking outward in this way, we find a rich vein of geometry,
topology and number theory associated with our theme of knots and unknots.
\bigbreak

\noindent {\bf Acknowledgements.} The first author thanks the
National Science Foundation for support of this research under NSF Grant DMS-0245588.
It gives both authors pleasure to acknowledge the hospitality of the  Mathematisches Forschungsinstitut Oberwolfach, the University of Illinois
at Chicago and the National Technical University of Athens, Greece, where much of this research was conducted. We particularly thank
Slavik Jablan for conversations and for helping us, with his computer program LinKnot, to find some key omissions in our initial enumerations.

\section{Culprits}

Classical knot theory is about the classication, up to isotopy, of embedded curves in three dimensional space. 
Two curves embedded in three-dimensional space are said to be {\it isotopic} if there is a continuous family of embeddings starting with
one curve and ending with the other curve. This definition of isotopy captures our intuitive notion of moving a length of rope flexibly about
in space without tearing it. Given a closed loop of rope in space, it can sometimes be deformed in this way to a simple flat circle. In
such a case, we say that the loop is {\it unknotted}. The perennial difficulty in the theory of knots is that unknots are not at all easy to 
recognise. View Figure 6. Here is illustrated an unknot, but most people find the unknotting difficult to visualize, although
it will present itself if one makes a rope model of the figure and shakes the model.
\bigbreak

One can make physical models of knots with rope. One can make mathematical models of knots by formalizing the diagrams that one naturally 
draws for them. For an example, see Figures 6 and 7. Knot theory based on the diagrams is called {\it combinatorial knot
theory}. The diagrams are made into mathematical structures by regarding them as planar graphs with four edges  incident to a
vertex and extra structure to indicate the crossing of one strand over the other at such a point.
\bigbreak
	 
Combinatorial knot theory got its start in the hands of Kurt Reidemeister \cite{Reidemeister} who discovered a set of moves on
planar diagrams that capture the topology of knots and links embedded in three dimensional space. Reidemeister proved
that the set of diagrammatic moves shown in Figure 1 generate isotopy of knots and links. That is, he showed that if we have two
knots or links in three dimensional space, then they are ambient isotopic if and only if corresponding
diagrams for them can be obtained, one from the other, by a sequence of moves of the types shown in Figure 1.

$$ \picill2.8inby2.9in(UK1) $$
\begin{center}
{\bf Figure 1 - The Reidemeister Moves}
\end{center}

In this section, we discuss the unknotting of knot diagrams using the Reidemeister moves. Reidemeister
proved that two knots or links are topologically equivalent if and only if diagrams for them can be related by a finite sequence
of  Reidemeister moves. The moves are performed on the diagrams in such a way that the rest of the diagram is left fixed. The
transformations shown in Figure 1 are performed locally. This means that when one searches for 
available Reidemeister moves on a diagram, one searches for one-sided, two-sided, and three-sided regions with the 
approriate crossing patterns as shown in Figure 1. It is important to note that knot and link diagrams are taken to be on the
surface of a two dimensional sphere that is standardly embedded in three dimensional space. Thus the outer region of a planar
diagram is handled just like any other region in that diagram, and can be the locus of any of the moves.
\bigbreak

Here is an example of a knot diagram (originally due to Ken Millett \cite{Ken}), in Figure 2. We
like to call this diagram the ``Culprit."  The Culprit is a knot diagram that  represents the unknot, but as a diagram, and 
using only the Reidemeister moves, it must be made more complicated before it can be simplified to an unknotted circle. We
measure the complexity of a knot or link diagram by the number of crossings in the diagram. Culprit has $10$ crossings, and in
order to be undone, we definitely have to increase the number of crossings before decreasing them to zero. The reader can verify
this for himself by checking each region in the diagram of the Culprit. A simplifying Reidemeister II move  can occur only on a
two-sided region, but no two-sided region in the diagram admits such a move. Similarly on the Culprit  diagram there are no
simplifying Redeimeister I moves and there are no Reidemeister III moves (note that a III move does not change the
complexity of the diagram). We view the diagram of the Culprit and other such examples as resting on the  surface of the
two-dimensional sphere. Thus the outer region of the diagram counts as much as any other region in this search for simplifying
moves.

$$ \picill1.5inby1.9in(UK2) $$
\begin{center}
{\bf Figure 2  - The Culprit }
\end{center}

$$ \picill1.5inby1.9in(UK3) $$
\begin{center}
{\bf Figure 3 - An Easy Unknot }
\end{center}

In contrast to the Culprit, view Figure 3 for an example of an easy unknot. This one simplifies after a Reidemeister
moves of type III unlock futher moves of type I and type II. The reader can probably just look at the diagram in Figure 3 and see
at once that it is unknotted.
\bigbreak 

View Figure 4 for an unknotting sequence for the Culprit. Notice that we undo it by swinging the arc that passes underneath
most of the diagram outward, and that in this process the number of crossings in the intermediate diagrams increases. In the
diagrams of Figure 4 the largest increase is to a diagram of $12$ crossings. This is the best possible result for this diagram.

$$ \picill8inby2.8in(UK4) $$
\begin{center}
{\bf Figure 4 - The Culprit Undone }
\end{center}

We shall call a diagram of the unknot {\em hard} if it has the following three properties:
\begin{enumerate}
\item There are no simplifying Type I moves on the diagram.
\item There are no simplifying Type II moves on the diagram.
\item There are no Type III moves on the diagram.
\end{enumerate}

Hard unknot diagrams have to be made more complex before they will simplify to the unknot, if we
use Reidemeister moves. It is an unsolved problem just how much complexity can be forced by a hard unknot.

\begin{defn} \rm
Define the {\it complexity} of a diagram $K$ to be the number of crossings,$C(K),$ of that diagram. Let $K$ be a hard unknot
diagram. Let $K'$ be a diagram Reidemeister equivalent to $K$ such that $K'$ can be simplified to the unknot.
For any unknotting sequence of Reidemeister moves for $K$ there will be a diagram $K'_{max}$ with a maximal number of
crossings. Let $Top(K)$ denote the minimum of $C(K'_{max})$ over all unknotting sequences for $K.$ Let
$$R(K) = Top(K)/C(K)$$
be called the {\em recalcitrance} of the hard unknot diagram $K.$ Very little is known about $R(K).$
\end{defn}

In the case of our Culprit, $Top(K) = 12$ while, $C(K) = 10.$ Thus $R(K) = 1.2.$ A knot that can 
be simplified with no extra complexity has recalcitrance equal to $1.$  We shall have more to say about the recalcitrance
in Section 13.
\bigbreak

In \cite{HassLagarias} Hass and Lagarias show that there exists a positive constant $c_1$ such that for each $n>1$ any unknotted
diagram $D$ with $n$ crossings can be transformed to the trivial knot diagram using less than or equal to $2^{c_{1}n}$ Reidemeister
moves.  As a corollary to this result, they conclude that any unknotted diagram $D$ can be transformed by Reidemeister
moves to the trivial knot diagram through a sequence of knot diagrams each of which has at most $2^{c_{1}n}$ crossings.
In our language, this result says that $$R(D) \le \frac{2^{c_{1}C(D)}}{C(D)}$$ where $C(D)$ is the number of crossings
in the diagram $D.$ The authors of \cite{HassLagarias} prove their result for $c_{1} = 10^{11},$ but remark that this is surely
too large a constant. Much more work needs to be done in this domain. We should also remark that the question 
of the knottedness of a knot is algorithmically decideable, due to the work of 
Haken and Hemion (See \cite{Hemion}). This algorithm is quite complex, but its methods are used in the work of Hass
and Lagarias.
\bigbreak

For work on unknots related to braids
see \cite{Morton}. See also the work of Dynnikov \cite{Dynnikov} for a remarkable diagramming system for knots and links
in which unknot diagrams can be detected by simplification only. Finally, another approach to detecting unknots is the use of
invariants of knots and links. For example, it is conjectured that the Jones polynomial \cite{Jones,KState} of a knot diagram
is equal to one only when that knot is unknotted. 
\bigbreak

In fact, we can now do considerably better in estimating $R(K).$ In \cite{HK} it is shown that
if the diagram $K$ represents the unknot and if $K$ is presented with a height function in the plane
(a Morse diagram) with $m(K)$ maxima and $c(K)$ crossings then one can unknot $K$ by Reidemeister moves using a diagram $K'$ with no more than $(c(K) + 2m(K) - 2)^{2}$ crossings.
This means that
$$R(K) \le \frac{(c(K) + 2m(K) - 2)^{2}}{c(K)}.$$
The problem of bounding the number
of Reidemeister moves needed to unknot $K$ is still difficult.
A technique of working with large-scale moves is formalized in the work of Dynnikov \cite{Dynnikov} where
he shows, using special grid diagrams and a complete set of moves distinct from Reidemeister moves, how to unknot diagrams  without making the special diagrams more complex. Dynnikov's
methods are used in proving the results in \cite{HK}.
\smallbreak

$$ \picill5.5inby2.9in(UK5) $$
\begin{center}
{\bf Figure 5 - The Culprit Analysis }
\end{center}

One purpose of this paper is to give infinite classes of hard unknots by employing an insight about the structure of our Culprit,
and generalizing this insight into results about the structure of tangles whose numerators are unknotted. In the course of 
this investigation, we shall obtain results about the collapse of tangle numerators to smaller knots and corresponding results
about the non-uniqueness of solutions to certain tangle equations. These results are of interest in working with the tangle model
of DNA recombination. See Section 8.
\bigbreak

In order to see the Culprit in a way that allows us to generalize him, we shall use the language and technique of the theory of 
tangles. The next sections describe a bit of basic tangle theory, but we shall now analyze the Culprit using this language, to 
illustrate our approach. The reader familiar with the language of tangles will have no difficulty here. {\it Other readers may wish
to  read the next section and then come back to this discussion.}
\bigbreak

View Figure 5. In Figure 5 we have a drawing $C$ of the Culprit and next to it, we have a drawing $C'$ of the result of part of
Figure 4, where the under-crossing arc in $C$ has been isotoped along (making the complexity rise) until it has been drawn outside
the rest of  the diagram. Concentrate on $C'.$ Notice that we can cut $C'$ into two pieces, as shown in Figure 5. These two
pieces, $A$ and 
$B$, are rational tangles (see the sections below) and this cutting process shows that $C'$ is the numerator closure of the 
tangle sum of $A$ and $B.$ This is written $C' = N(A + B).$ Each rational tangle $T$ has a rational fraction $F(T)$ that tells
all about it. In this case, $$F(A) = \frac{1}{-1 + \frac{1}{-3}} = -3/4$$ and $$F(B) = \frac{1}{1 + \frac{1}{2}} = 2/3.$$
We know that the numerator of $A + B$ is unknotted and we would like to understand why it is unknotted. We notice that 
the sum of the fractions of $A$ and $B$ is $-3/4 + 2/3 = -1/12.$ Thus the numerator of the sum of the fractions of 
$A$ and $B$ is $-1.$ Does this $-1$ imply the unknottedness of the numerator of $A + B?$ Well, the answer is that it does, and
that will be the subject of much of the rest of this paper. See our Theorem 5 in Section 3.
\bigbreak

The next example is illustrated in Figures 6 and Figure 7. Figure 6 is a three-dimensional rendering of the unknot
diagram shown in  Figure 7. This example is interesting psychologically because it looks knotted to most observers. The
three dimensional picture is a frame from a deformation of this example in the energy minimization program KnotPlot
\cite{Scharein}. It is illuminating to run the program on this example and watch the knot self-repel and undo itself. 
A little tangle arithmetic shows that this knot is unknotted. The diagram shown in Figure 7 is the
numerator of a sum of two tangles (the dotted line in this figure indicates the tangle decomposition). One tangle is $[1/2].$
The other is equivalent to the sum
$[1] + [-1] + [-1/3] = [-1/3].$ Since 
$1/2  - 1/3 = 1/6,$ this knot is unknotted by our Theorem 5. Note that this unknot example is not hard. It has type III moves
available on the diagram.

$$ \picill2inby2in(UK6.EPSF) $$
\begin{center}
{\bf Figure 6 - Three Dimensional Unknot}
\end{center}

$$ \picill2inby2in(UK7) $$
\begin{center}
{\bf Figure 7 - Three Dimensional Unknot Diagram}
\end{center}


Finally, we give an example of a hard unknot that is of a different type than the sort that we are considering in this paper.
In Figure 8 we show a hard unknot whose diagram is not the closure of the sum of two rational tangles.  There are many knots
and many hard unknots. We only scratch the surface of this subject.

$$ \picill3inby2.5in(UK8.EPSF) $$
\begin{center}
{\bf Figure 8 - Another Hard Unknot}
\end{center}

\section{Rational Tangles, Rational Knots and Continued Fractions}

In this section we recall the subject of  rational tangles and 
rational knots and their relationship with the theory of continued fractions. By the term ``knots" we
will refer to both knots and links, and whenever we really mean ``knot" we shall emphasize it. Rational knots and links
comprise the simplest class of links. They are also known in the literature as Viergeflechte, four-plats  or
$2$-bridge knots depending on their geometric representation. The notion of a tangle was introduced in 1967 by
Conway \cite{C1} in his work on enumerating and classifying knots and links. 
\bigbreak

A {\it $2$-tangle} is a proper embedding of two unoriented arcs and a finite number of circles in a $3$-ball
$B^3,$ so that the four endpoints lie in the boundary of $B^3$. A {\it tangle diagram} is a regular projection
of the tangle on an equatorial disc of $B^{3}.$  By ``tangle" we will mean ``tangle
diagram".  A {\it rational tangle} is a special case of a $2$-tangle obtained by applying consecutive twists
on neighbouring endpoints of two trivial arcs. Such a pair of arcs comprise the $[0]$ or $[\infty]$ tangles 
(see Figure 10), depending on their position in the plane. We shall say that the rational tangle is in {\it
twist form} when it is obtained by such successive twists. For examples see  Figure 14. Conway defined the rational knots as
``numerator" or ``denominator" closures  of the rational tangles. See Figure 9.  Conway \cite{C1} also defined  {\it the
fraction} of a rational tangle to be a rational number or
$\infty,$ obtained via a continued fraction that is associated with the tangle. We discuss this construction below.

$$ \picill8inby1.7in(UK9)  $$
\begin{center}
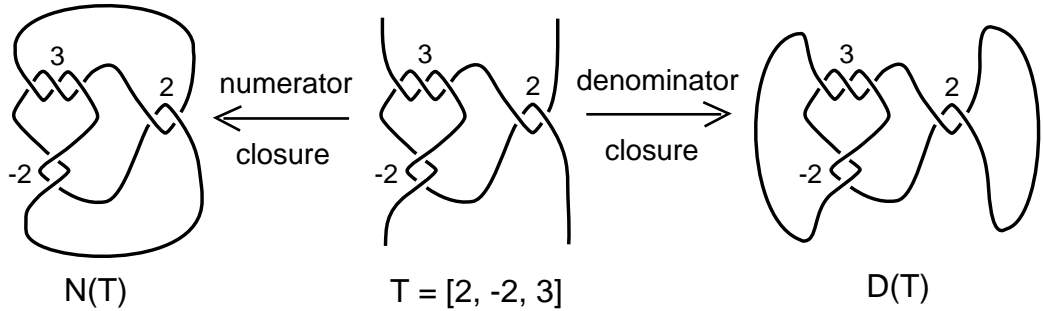
{\bf Figure 9 - A rational tangle and its closures to rational knots } 
\end{center} 

We  are interested in tangles up to isotopy. Two
tangles, $T, S$, in $B^3$ are {\it isotopic}, denoted by $T \sim S$,  if and only if any two  diagrams of
them have identical configurations of their four endpoints on the boundary of the
projection disc, and they differ by a finite sequence of the Reidemeister moves
\cite{Reidemeister,Rd2}, which take place in the interior of the  disc.  Of course, each  twisting
operation used in the definition of a rational tangle changes the isotopy class of  the tangle to which
it is applied. Rational tangles are classified by their fractions by
means of the following theorem, different proofs of which are given in \cite{Mo}, \cite{Bu,BZ}, \cite{GK2}
and \cite{KL1}.   

\begin{th}[Conway, 1975]\label{Conway}{ \ Two rational tangles are isotopic if and
only if they have the same fraction.
  } \end{th}
 
More than one rational tangle can yield the same or isotopic  
rational knots, and the equivalence relation between the rational tangles is mapped  
into an arithmetic equivalence of their corresponding fractions. Indeed we have:

\begin{th}[Schubert, 1956]\label{Schubert1}{ \ Suppose that rational tangles with
fractions $\frac{p}{q}$ and $\frac{p'}{q'}$ are given ($p$ and $q$ are relatively prime.
Similarly for $p'$ and $q'$.) If  $K(\frac{p}{q})$ and $K(\frac{p'}{q'})$
denote the corresponding rational knots obtained by taking numerator closures of
these tangles, then $K(\frac{p}{q})$ and $K(\frac{p'}{q'})$ are isotopic 
if and only if

\begin{enumerate}
\item $p=p'$ and
\item either $q\equiv q'\, mod \, p$  \ or \ $qq'\equiv 1\, mod \, p.$
\end{enumerate}
} \end{th}

Different proofs of Theorem 2 are given in \cite{Sch2}, \cite{Bu}, \cite{KL2}.

\subsection{Rational Tangles and their Invariant Fractions}

We shall now recall from \cite{KL1} the main properties of rational tangles and of continued fractions, which
illuminate the classification of rational tangles.  
The elementary rational tangles are displayed as either
horizontal or vertical twists, and they are enumerated by integers or their inverses, see Figure 10. 

\smallbreak

The crossing types of
$2$-tangles (and of unoriented knots) follow the checkerboard  rule: shade the regions of the tangle in two
colors, starting from the left outside region with grey, and so that
adjacent regions have different colors. Crossings in the tangle
are said to be of ``positive type" if they are arranged with respect to the shading 
as exemplified in Figure 10 by the tangle $[+1],$ i.e. they have the  region on the right
shaded as one walks towards the crossing along the over-arc. Crossings of the reverse
type are said to be of ``negative type" and they are exemplified in Figure 10 by the
tangle $[-1].$ The reader should note that our crossing type  conventions are
the opposite of those of Conway in \cite{C1}. Our
conventions agree with those of  Ernst and Sumners \cite{Sumners1,Sumners2}, \cite{Su} which in turn follow the
standard conventions of  biologists.

$$ \picill4.7inby1.7in(UK10) $$
\begin{center}
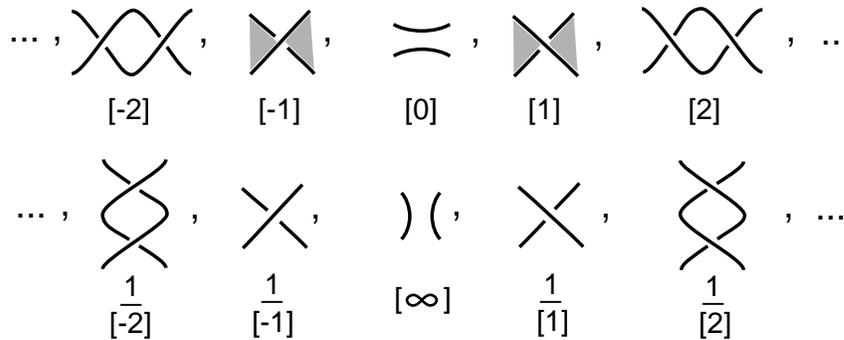
{\bf Figure 10 - The elementary rational tangles and the types of crossings } 
\end{center}

In the class of $2$-tangles we have the non-commutative operations 
{\it addition} and {\it multiplication},  as illustrated in Figure 11, which are denoted by
``$+$" and  ``$*$" respectively. These operations are well-defined up to isotopy. A rational
tangle in twist form is created inductively by consecutive additions of the tangles
$[\pm 1]$ on the right or on the left and multiplications  by the  tangles $[\pm 1]$ at
the bottom or at the top, starting from the tangle $[0]$ or $[\infty]$.  Since the very
first  crossing  can be equally seen as horizontal or  vertical, we may always assume 
that we start twisting from the tangle $[0]$. In order to read out a rational tangle  we transcribe
it as an algebraic sum using  horizontal and vertical twists. For example, Figure 9 illustrates the tangle 
$(([3]*\frac{1}{[-2]})+[2])$, see top of Figure 12, while Figure 14 illustrates a twist form of the same
tangle: $[1]+([1]*[3]*\frac{1}{[-3]})+[1]$. 

\smallbreak

Note that addition and multiplication do not, in general, preserve the class of rational tangles. For
example, the $2$-tangle $\frac{1}{[3]} + \frac{1}{[3]}$ is not rational. The sum (product) of two rational
tangles is rational if and only if one of the two consists in a number of horizontal (vertical) twists. 
  
\bigbreak

$$\vbox{\picill4.6inby1.6in(UK11)  }$$
\begin{center}
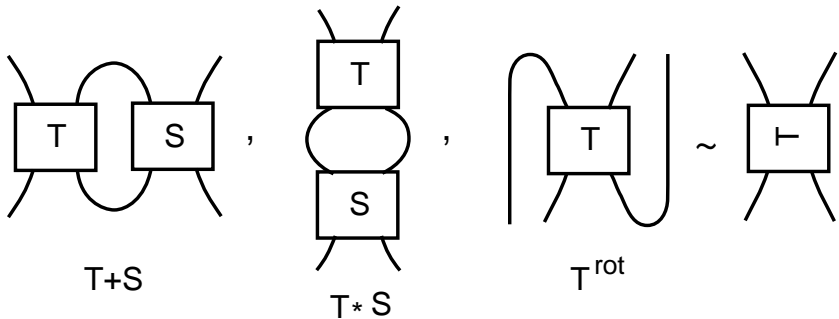
{\bf Figure 11 - Addition, multiplication and rotation of $2$-tangles } 
\end{center}

The {\it mirror image} of a tangle $T,$  denoted $-T,$ is $T$ with all crossings switched. For
example, $-[n] = [-n]$ and  $ -\frac{1}{[n]} =\frac{1}{[-n]}.$  Then, the {\it subtraction} is defined as
$T-S := T + (-S)$.  The {\it rotation} of $T$, denoted 
$T^{rot}$, is obtained by rotating $T$ on its plane counterclockwise by $90^0.$ The {\it inverse} of $T$
is defined to be  ${-T}^{rot}$. 
Thus, inversion is accomplished by rotation and mirror image. Note that $T^{rot}$ and the inverse of $T$
are in general not isotopic to $T$ and they are order 4 operations. But for rational  tangles the inversion
is an operation of order 2 (this follows from the flipping lemma discussed below). For this reason we shall denote the
inverse of a rational tangle
$T$ by
$1/T,$ and hence the rotation of the tangle $T$ will be denoted by $-1/T.$ This explains the
notation for the tangles $\frac{1}{[n]}.$

$$ \picill5inby6.8in(UK12) $$
\begin{center}
{\bf Figure 12 - Finding the Fraction } 
\end{center}

As we said earlier, there is a fraction associated to a rational tangle $R$ which characterizes its isotopy
class (Theorem 1). In fact, the fraction is defined for any $2$-tangle and always has the following three 
properties. These suffice for computing the fraction $F(R)$ inductively for rational tangles:

\begin{enumerate}
\item $F([\pm1]) = \pm1.$
\item $F(T + S) = F(T) + F(S).$
\item $F(T^{rot}) = -1/F(T).$
\end{enumerate}

In Figure 12 we illustrate this process by using only these three rules to compute a specific tangle fraction.
In the following discussion we discuss the fraction in more detail and how it is related to the continued fraction
structure of the rational tangles.
\bigbreak

For rational tangles in twist form two types of isotopy moves suffice for their
study: the {\it  flypes} and the {\it  swing moves},  illustrated in Figures 13 and Figure 17 below. Both
moves apply on
$2$-subtangles, fixing their endpoints.

$$ \picill3.65inby.6in(UK13) $$

\begin{center}
{\bf Figure 13 - The flype moves } 
\end{center}

Starting from a rational tangle in twist form and using obvious flypes on appropriate subtangles one can
always bring the twists all to the right and to the bottom of the tangle. See Figure 14 for an example.  
We shall then say that the  rational tangle is in {\it standard form}. 

For the study of rational knots it is easier to use another way of representing an abstract rational tangle in
standard form, illustrated in Figure 14. This is the {\it $3$-strand-braid representation}. As illustrated in Figure 14,
the $3$-strand-braid representation is obtained from the standard representation by planar rotations of
the vertical sets of crossings, thus creating a lower row of horizontal crossings. Note that the
type of crossings does not change by this planar rotation. Indeed the checkerboard coloring convention for the
crossing signs identifies the signs as unchanged. Nevertheless, the crossings on the lower row of the braid representation
appear to be of opposite sign, since when we rotate them to the vertical position we obtain crossings of
the opposite type in the local tangles. We shall use both ways of representation for extracting the
properties of rational knots and tangles.

$$ \picill8inby1.8in(UK14)  $$

\begin{center}
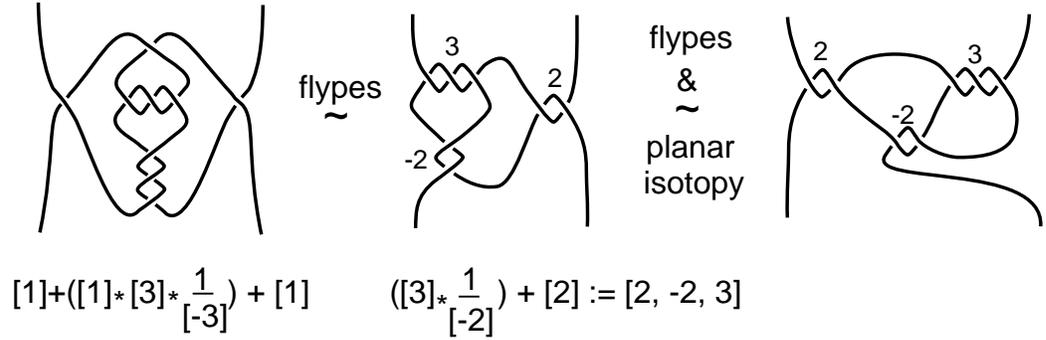
{\bf Figure 14 - A rational tangle in twist form converted to its standard form and to its $3$-strand-braid representation} 
\end{center} 

One can associate to a rational tangle in standard form a vector of integers $(a_1, a_2, \ldots, a_n), $
where the first entry denotes the place where the tangle starts untwisting and the last
entry where it  begins to twist.  For example the tangle of Figure 9 corresponds  to
the vector $(2, -2, 3)$.

\smallbreak 

Note that the set of twists of a rational tangle may be always assumed  {\it odd}. Indeed,  let  $n$ be even
and let the left-most twist $[a_1]$ be on the upper part of the braid representation. 
Then, the right-most crossing of the last twist $[a_n]$
may be assumed upper, so that 
$[a_n]$ can break into $a_n-1$ lower crossings and one upper. Up to the ambiguity of the
 right-most crossing, the vector associated to a rational tangle is {\it unique}, i.e. $(a_1, a_2,
\ldots, a_n) =  (a_1, a_2, \ldots, a_n -1, 1),$ if $a_n >0,$ and $(a_1, a_2, \ldots, a_n) =  (a_1, a_2,
\ldots, a_n +1, -1),$ if $a_n <0$. See Figure 15.

$$ \picill4.5inby.9in(UK15) $$

\begin{center}
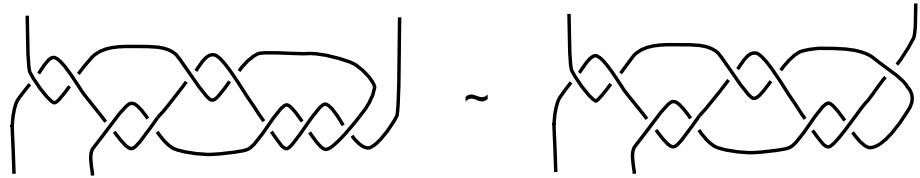
{\bf Figure 15 - The ambiguity of the first crossing  } 
\end{center} 

\bigbreak

Another move that can be applied to a $2$-tangle  is a {\it flip}, its rotation in space  by $180^0$. We
denote $T^{hflip}$ a horizontal flip (rotation around a horizontal axis on the plane of $T$) and
$T^{vflip}$ a vertical flip. See Figure 16 for illustrations. Note that a flip switches the endpoints of the
tangle and, in general, a flipped tangle is not isotopic to the original one. {\it Rational tangles have the
remarkable property that they are isotopic to their horizontal or vertical flips.} We shall refer to this as
the {\it  Flipping Lemma}.

$$ \picill3.5inby2.4in(UK16) $$

\begin{center}
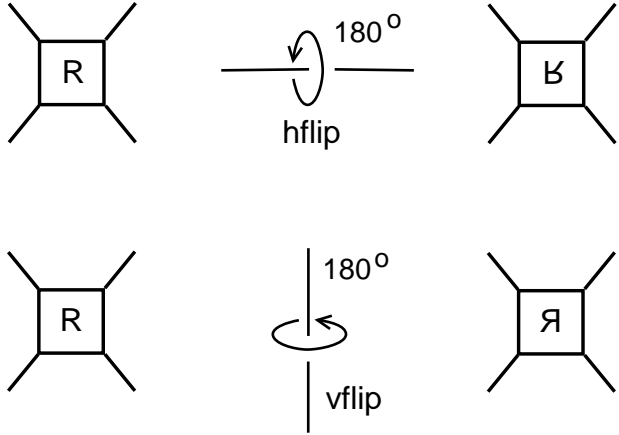
{\bf Figure 16 - The horizontal and the vertical flip } 
\end{center}

A consequence of the Flipping Lemma is that addition and multiplication by  $[\pm
1]$ are commutative.
Another consequence of the Flipping Lemma is that rotation and inversion of rational tangles each have order 
two. In particular, rotation is defined via a ninety degree turn of the tangle either to the left or to the
right.   With this in mind the reader  can easily deduce the formula below:  
$$ 
T* \frac{1}{[n]} = \frac{1}{[n] + \frac{1}{T}}
$$ 
Indeed, rotate $ T* \frac{1}{[n]}$ by ninety degrees and note that it becomes $-[n] - \frac{1}{T}$.  Use
this to deduce that the original tangle is the negative reciprocal of this tangle.  This formula implies
that the two operations:
 addition of $[+1]$ or $[-1]$  and inversion  between rational
tangles suffice for generating the whole class of rational tangles.
 As for the fraction, we have the corresponding formula
$$ 
F(T* \frac{1}{[n]}) = \frac{1}{n + \frac{1}{F(T)}}.
$$ 
 The above equation for tangles leads to 
the fact that  a rational tangle in standard form can be described
 algebraically by a continued fraction built from the integer tangles  $[a_1],$
 $[a_2],$ $\ldots, [a_n]$  with all numerators equal to~$1,$ namely by an expression of
 the type:

 \[ [[a_1],[a_2],\ldots ,[a_n]] \, := \, [a_1]+ \frac{1}{[a_2]+\cdots +
 \frac{1}{[a_{n-1}]
 +\frac{1}{[a_n]}}} \]

 \noindent for $a_2, \ldots, a_n \in \ZZ - \{0\}$ and $n$ even or odd. We allow  $[a_1]$  to be the
tangle $[0].$  Then, a rational tangle  is said to be in {\it continued fraction form}. 
\bigbreak

\noindent {\it We shall abbreviate the expression
$[[a_1],[a_2],\ldots ,[a_n]]$ by writing $[a_1,a_2,\ldots ,a_n],$} and later will use the latter expression for a numerical continued 
fraction as well. There should be no ambiguity between the tangle and numerical interpretations, as these will be clear from context. 
Figure 9 illustrates the rational tangle $[2,-2,3]$.   
\bigbreak

\noindent Let  $T=[a_1,a_2,\ldots,a_n].$  
The following statements are now straightforward. 
 \vspace{.15in}

 \noindent $\begin{array}{lrcl}

  1. \, \, \, \, \, \,  T + [\pm 1]   =  [a_1 \pm 1,a_2,\ldots,a_n],   \\

  2. \, \, \, \, \, \,  \frac{1}{T}   =  [0,a_1,a_2,\ldots,a_n],   \\ 

  3. \, \, \, \, \, \,   -T   =   [-a_1,-a_2,\ldots,-a_n]. \\

  4. \, \, \, \, \, \,   T   =   [a_1,a_2,\ldots,a_n -1,1], \ \ \mbox{if $a_n >0,$ and}  \\

     \, \, \, \, \, \, \, \, \, \,   T  =  [a_1,a_2,\ldots,a_n +1,-1],  \ \ \mbox{if $a_n <0.$} \\

  5. \, \, \, \, \, \,  [a_1,\ldots, a_{i-1}, b_i, 0, c_i, a_{i+1},\ldots, a_n] = [a_1,\ldots, a_{i-1}, b_i + c_i, a_{i+1}, \ldots, a_n].\\

\end{array}$
\bigbreak

From the above discussion it makes sense to assign to a rational tangle in standard
form, $T=[[a_1],[a_2],\ldots,[a_n]],$ for 
$ a_1 \in \ZZ, \ a_2, \ldots, a_n \in \ZZ - \{0\}$ and $n$
 even or odd, the numerical continued fraction 
 \[ F(T) = F([[a_1],[a_2],\ldots,[a_n]]) = [a_1, a_2, \ldots, a_n] := a_1+ \frac{1}{a_2+\cdots +
\frac{1}{a_{n-1}
 +\frac{1}{a_n}}},  \] 

\noindent if $T \neq [\infty],$ and to assign $F([\infty]) := \infty = \frac{1}{0},$
as a formal expression. This rational number or infinity shall be called {\it the
fraction  of $T$}. With this definition of the fraction of a rational tangle, all of the above properties
of the continued fraction form of the tangle become properties of the fraction. We leave this translation
process to the reader.
\bigbreak

The subject of continued  fractions is of perennial interest to mathematicians. See for
example \cite{Kh},  \cite{O}, \cite{Ko}, \cite{W}. Here we  only consider continued fractions of the above
type, with all numerators equal to 1. As in the case of rational tangles we allow the term $a_1$  to be
zero. Clearly, the two simple algebraic operations {\it addition  of $+1$ or $-1$} and {\it inversion}
generate inductively the whole class of continued fractions, starting from
zero.  
\bigbreak

If a rational tangle $T$ changes by an isotopy, the associated continued fraction form may also
change.  However, the fraction is a topological invariant of $T$ and does not change. For example,
$[2,-2,3] = [1,2,2] = \frac{7}{5}$, see Figure 12. The fraction characterizes the isotopy class of $T$.  For
the isotopy type of a rational tangle
$T$ with fraction
$\frac{p}{q}$ we shall use the notation $[\frac{p}{q}]$. We have omitted here the proof of the invariance of
the fraction. The interested reader can consult \cite{C1}, \cite{GK2}, \cite{KL1} for various proofs of this
fact.
\bigbreak

The key to the exact correspondence of fractions and rational tangles lies in the construction of a canonical
alternating form  for the rational tangle. This is obtained as follows:
A tangle is
said to be {\it alternating} if it is isotopic to a tangle where the crossings  alternate from under to over as we go along any
component or arc of the weave. Similarly, a knot is alternating if it possesses an alternating diagram.
According to the checkerboard shading, the weave can alternate only if any two adjacent crossings are of the
same type, and this propagates to the whole diagram. As we mentioned
earlier, the second type of isotopy moves that suffice for studying rational tangles are the swing moves,
illustrated on an example in Figure 17. Using the swing moves and an induction argument, it is easy to show
that rational tangles (and rational knots) are alternating.

\smallbreak
 We shall say that the  rational tangle
$S=[\beta_1, \beta_2, \ldots, \beta_m]$ is in  {\it  canonical form} if $S$ is alternating and $m$ is odd.
From the above, $S$ alternating implies that the $\beta_i$'s are all of the same sign. It turns out that the
canonical form for $S$ is unique. In Figure 17 we bring our working rational tangle $T=[2,-2,3]$ to its
canonical form $S=[1,2,2]$. As noted above, $F(T)=F(S)=\frac{7}{5}$.

$$ \picill8inby1.2in(UK17) $$

\begin{center}
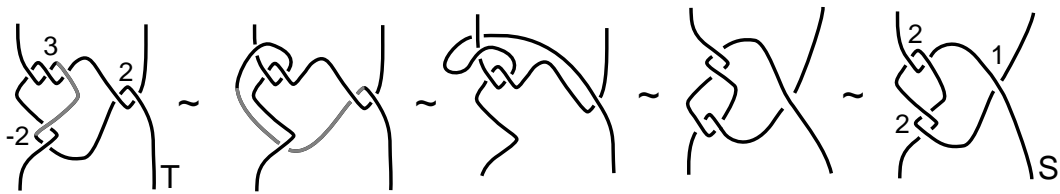
{\bf Figure 17 - Reducing to alternating form using the swing moves} 
\end{center} 

On the other hand, by Euclid's algorithm and keeping all remainders of the same sign, one can show that
 every  continued fraction $[a_1, a_2, \ldots, a_n]$
can be transformed to a unique canonical form $ [\beta_1, \beta_2, \ldots, \beta_m],$ 
where all $\beta_i$'s are positive or all negative integers and $m$ is odd. For example,
$[2,-2]=[1,1,1]=\frac{3}{2}$. There is also an algorithm that can be applied directly to the initial
continued fraction to obtain its canonical form, which works in parallel with the algorithm for the
canonical form of rational tangles. Indeed, we have:  

\begin{prop}
The following identity is true for continued fractions and it is also a topological
equivalence of the corresponding tangles: 
$$[..., a,-b,c,d,e,...] = [..., (a - 1) , 1, (b-1), -c, -d, -e, ...].$$
This identity gives a specific inductive procedure for reducing a continued fraction to all positive or all 
negative terms. In the case of transforming to all negative terms, we can first flip all signs and work with the mirror image.
Note also that $$[..., a,b,0,c,d,e, ...] = [..., a,b + c,d,e, ...]$$ will be used in these reductions.
\end{prop}

\noindent {\bf Proof.} The technique for the reduction is based on the formula
$$a + 1/(-b) = (a-1) + 1/(1 + 1/(b-1)).$$
If $a$ and $b$ are positive, this formula allows the reduction of negative terms in a continued fraction.
The identity in the Proposition follows immediately from this formula.  
$\hfill \Box$

\subsection{Rational Knots and Continued Fractions}

By joining with simple arcs the two upper and the two lower endpoints of a $2$-tangle $T,$ we obtain a 
knot called the {\it Numerator} of $T$, denoted by $N(T)$. A rational knot is defined to be the
numerator of a rational tangle. Joining with simple arcs  each pair of the corresponding top and bottom
endpoints of $T$ we obtain the {\it Denominator} of $T,$  denoted by $D(T)$, see Figure 9. We have $N(T) =
D(T^{rot})$ and $D(T) = N(T^{rot}).$  As we shall see in the next section,  the numerator closure of the sum of
two rational tangles  is still a rational knot. But the denominator closure of the sum of
two rational tangles  is not necessarily a rational knot, think for example of the sum 
$\frac{1}{[3]} + \frac{1}{[3]}.$

\bigbreak

Given two different rational tangle types  $[\frac{p}{q}]$ and $[\frac{p'}{q'}]$, when do
they close to isotopic rational knots? The answer is given in Theorem 2. Schubert  classified
rational knots by finding canonical forms via representing them as $2$-bridge knots.  In \cite{KL2} we give a
new combinatorial proof of Theorem 2, by posing the question: given a rational knot diagram, at which places may one cut it so 
that it opens to a rational tangle?
We then pinpoint two distinct categories of cuts that
represent the two cases of the arithmetic equivalence of Schubert's theorem. The first case corresponds to
the {\it  special cut}, as illustrated in Figure 18.
 The two tangles $T  = [-3] $ and $S = [1] + \frac{1}{[2]} $ are non-isotopic by the Conway
 Theorem, since $F(T) =  -3 = 3/-1, $ while $F(S) =1 + 1/2 = 3/2.$ But they have isotopic
 numerators:  $N(T) \sim N(S),$ the left-handed trefoil. Now  $-1 \equiv 2 \, mod \, 3,$  confirming
 Theorem \ref{Schubert1}. See \cite{KL2} for
a complete analysis of the special cut.

$$ \picill3.45inby1.6in(UK18) $$

\begin{center}
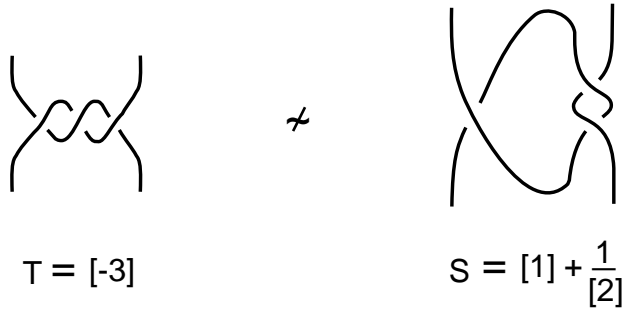
{\bf Figure 18 - An example of the special cut } 
\end{center}

The second case of Schubert's equivalence corresponds to the {\it  palindrome cut}, an example of which is
illustrated in Figure 19. Here we see that the tangles 

$$ T=[2,3,4] = [2] + \frac{1}{[3] + \frac{1}{[4]}} $$

\noindent and 
$$ S=[4,3,2] = [4] + \frac{1}{[3] + \frac{1}{[2]}} $$

\noindent both have the same numerator closure. Their corresponding fractions are 
 $$ F(T)= 2 + \frac{1}{3 + \frac{1}{4}} = \frac{30}{13} \mbox{ \ \ and \ \  } F(S) = 4 +
\frac{1}{3 + \frac{1}{2}} = \frac{30}{7}.$$
 Note that \  $7\cdot 13\equiv 1 \, mod \, 30. $

$$ \picill5inby2.2in(UK19) $$

\begin{center}
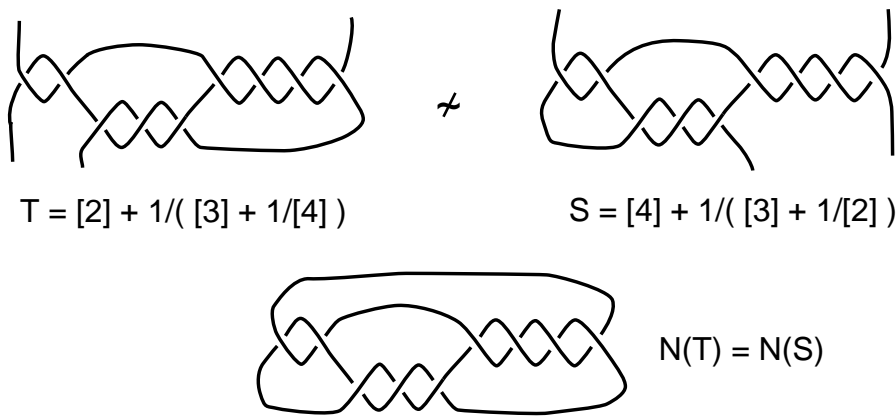
{\bf Figure 19 - An instance of the palindrome equivalence  } 
\end{center} 

In the general case if 
$ T=[a_1,a_2, \ldots, a_n], $ we shall call the tangle $ S=[a_n,a_{n-1}, \ldots,
a_1]$ {\it the  palindrome of $T$.} Clearly these tangles have the same numerator. 
In order to check the arithmetic in the general case of the palindrome cut  we need to 
generalize this pattern to arbitrary continued fractions and their palindromes (obtained
by reversing the  order of the terms).  
\bigbreak

The next Theorem is a known result about continued fractions. See \cite{KL1},
\cite{Sie} or \cite{Kaw}.  We shall give here our proof of this
statement. For this we will first present a way of evaluating continued fractions via
$2\times 2$ matrices (compare with \cite{Fr}, \cite{Ko}). This method of
evaluation is crucially important for the rest of the paper. We define matrices
$M(a)$ by the formula

$$M(a) =
\left( \begin{array}{cc}
a & 1 \\
1 & 0\\
\end{array} \right).$$

\noindent  These matrices $M(a)$ are said to be the {\it generating matrices} for
continued  fractions,  as  we have: 

\begin{th}
[The matrix product interpretation for continued fractions]
\bigbreak

Let $\{ a_{1},a_{2},\ldots ,a_{n} \}$
be a collection  of $n$ integers, and let 
$$ \frac{P}{Q} = [a_{1},a_{2},\ldots,a_{n}]$$ and
$$ \frac{P'}{Q'} =[a_{n},a_{n-1},\ldots ,a_{1}].$$ Then 
$P=P'$ and  $QQ' \equiv  (-1)^{n+1}\, mod \, P.$
\smallbreak

In fact, for any sequence of 
integers $\{  a_1,a_2,\ldots ,a_n \}$ 
the value of the corresponding continued fraction
$$\frac{P}{Q} = [a_{1},a_{2},\ldots ,a_{n}]$$
is given through the following  matrix
product $$M = M(a_{1})M(a_{2})\cdot\cdot\cdot M(a_{n})$$ via the identity
$$M =
\left( \begin{array}{cc}
P & Q' \\
Q & U \\
\end{array} \right)$$ where this matrix also gives the 
evaluation of of the palindrome continued fraction  
$$[a_{n},a_{n-1},\ldots ,a_{1}] = \frac{P}{Q'}.$$  
\end{th}

\noindent {\em Proof.}  
Let $$\frac{R}{S} = [a_{2},a_{3},\ldots ,a_{n}].$$ Then
$$\frac{P}{Q} = [a_{1},a_{2},\ldots ,a_{n}] = a_{1} + \frac{1}{\frac{R}{S}} = a_{1} + \frac{S}{R}$$
$$= \frac{Ra_{1} + S}{R}.$$
By induction we may assume that 
$$M(a_{2})M(a_{3}) \cdots M(a_{n}) =
\left( \begin{array}{cc}
R & S' \\
S & U \\
\end{array} \right).$$ 
Hence
$$M(a_{1})M(a_{2})\cdot\cdot\cdot M(a_{n}) = \left( \begin{array}{cc}
a_{1} & 1 \\
1 & 0\\
\end{array} \right)
\left( \begin{array}{cc}
R & S' \\
S & U \\
\end{array} \right) =
\left( \begin{array}{cc}
a_{1}R + S & a_{1}S' + V \\
R & S' \\
\end{array} \right).$$
This proves by induction that 
$$M(a_{1})M(a_{2})\cdot\cdot\cdot M(a_{n}) = 
\left( \begin{array}{cc}
P & Q' \\
Q & U \\
\end{array} \right)$$
where $$\frac{P}{Q} = [a_{1},a_{2},\ldots ,a_{n}].$$

To see the result about the palindrome continued fraction, just note that
if $M^{T}$ denotes the transpose of a square matrix $M,$ then
with $$M = M(a_{1})M(a_{2})\cdot\cdot\cdot M(a_{n}),$$ 
$$M^{T} = M(a_{n})M(a_{n-1}) \cdots M(a_{1}),$$
from which the statement about the palindrome follows.
Note also that $Det(M) = (-1)^{n}$ since $M$ is a product of $n$ matrices of determinant equal to minus 
one,
and $Det(M) = PU - QQ'.$ Hence $$PU - QQ' = (-1)^{n}.$$ This last equation implies the congruence stated
 at the beginning of the
Theorem, and completes the proof of the Theorem.
$\hfill \Box$

\section{Sums of Two Rational Tangles}

In this section we prove that the numerator of the sum of two rational tangles is a rational knot or link. We characterize 
the knot or link that emerges from this process.

$$ \picill5inby3.1in(UK20) $$
\begin{center}
{\bf Figure 20 - The numerator of a sum of rational tangles is a rational link}
\end{center}

\begin{th}
[Addition of Rational Tangles]\label{ratlknot} 
Let $\{ a_{1},a_{2},\ldots ,a_{n} \}$
be a collection  of integers, so that  $$ \frac{P}{Q} = [a_{1},a_{2},\ldots
,a_{n}].$$ Let $\{ b_{1},b_{2},\ldots ,b_{m} \}$
be another collection  of integers, so that  $$ \frac{R}{S} = [b_{1},b_{2},\ldots
,b_{m}].$$ Let $A = [\frac{P}{Q}]$ and $B = [\frac{R}{S}]$ be the corresponding
rational tangles. Then the knot or link $N(A + B)$ is rational, and in fact
$$N(A + B) = N([a_{n}, a_{n-1}, \ldots , a_{2}, a_{1} + b_{1}, b_{2}, \ldots , b_{m}]).$$
\end{th}

\noindent {\em Proof.}  View Figure 20. In this figure we illustrate a special case of the Theorem.
The geometry of reconnection in the general case should be clear from this illustration.
$\hfill \Box$
\bigbreak

The next result tells us when we get the unknot.
\bigbreak

\begin{defn} \rm
Given continued fractions $\frac{P}{Q} = [a_{1},\ldots,a_{n}]$ and $\frac{R}{S} = [b_{1},\ldots,b_{m}]$, let  
$$[a_{1},\ldots,a_{n}] \, \sharp \, [b_{1},\ldots,b_{m}] = [a_{n},\ldots,a_{2}, a_{1} + b_{1}, b_{2},\ldots, b_{m}].$$
If $$\frac{F}{G} = [a_{n},\ldots,a_{2}, a_{1} + b_{1}, b_{2},\ldots, b_{m}],$$ we shall write
$$\frac{P}{Q} \, \sharp \, \frac{R}{S} = \frac{F}{G}.$$ Note that 
$\frac{F}{G}$ is a fraction  such that $N([\frac{F}{G}]) = N([\frac{P}{Q}] + [\frac{R}{S}]).$ 
\end{defn}

\begin{th}
Let $$ \frac{P}{Q} = [a_{1},a_{2},\ldots
,a_{n}]$$ and $$ \frac{R}{S} = [b_{1},b_{2},\ldots
,b_{m}]$$  be as in the previous Theorem. Then $$N([\frac{P}{Q}] + [\frac{R}{S}])$$ is unknotted
if and only if $PS + QR = \pm 1.$
\end{th}

\noindent{\em Proof.} Let $M(\vec{a})$ denote the product of matrices
$$M(\vec{a}) = M(a_{1}) \cdots M(a_{n}).$$
We know from Theorem 3 that $$M(\vec{a}) = M(a_{1}) \cdots M(a_{n}) = 
\left( \begin{array}{cc}
P & Q' \\
Q & U \\
\end{array} \right),$$ and that 
$$M(\vec{b}) = M(b_{1}) \cdots M(b_{m}) =
\left( \begin{array}{cc}
R & S' \\
S & V \\
\end{array} \right).$$
Let $$F/G = [a_{n}, \ldots , a_{2}, a_{1} + b_{1}, b_{2}, \ldots ,b_{m}] = \frac{P}{Q}\sharp\frac{R}{S},$$
and let $$M(\vec{c}) = M(a_{n}) \cdots M(a_{2})M(a_{1} + b_{1})M(b_{2}) \cdots M(b_{m}).$$
Then we have, by Theorem 4, that 
$$N([\frac{P}{Q}] + [\frac{R}{S}]) = N([\frac{F}{G}]),$$ and 
$$M(\vec{c}) = M(a_{n}) \cdots M(a_{2}) M(a_{1} + b_{1}) M(b_{2}) \cdots M(b_{m}) =
\left( \begin{array}{cc}
F & G' \\
G & W \\
\end{array} \right).$$
Now note the identity
$$\left( \begin{array}{cc}
a_{1} & 1 \\
1 & 0 \\
\end{array} \right)
\left( \begin{array}{cc}
0 & 1 \\
1 & 0 \\
\end{array} \right)
\left( \begin{array}{cc}
b_{1} & 1 \\
1 & 0 \\
\end{array} \right) =
\left( \begin{array}{cc}
a_{1} + b_{1} & 1 \\
1 & 0 \\
\end{array} \right).$$
Thus
$$M(\vec{c}) = M(\vec{a})^T
\left( \begin{array}{cc}
0 & 1 \\
1 & 0 \\
\end{array} \right)
M(\vec{b}) =
M(\vec{a})^T M(\vec{b})'$$
where $M'$ denotes the matrix obtained from the 
$2 \times 2$ matrix $M$ by interchanging its two rows.
In particular, this formula implies that 
$$\left( \begin{array}{cc}
F & G' \\
G & W \\
\end{array} \right) =
\left( \begin{array}{cc}
P &  Q \\
Q' & U \\
\end{array} \right)
\left( \begin{array}{cc}
S & V \\
R & S' \\
\end{array} \right).$$
Hence $F = PS + QR.$
From this it follows that $N([\frac{F}{G}])$ is unknotted if and only if
$PS + QR = \pm 1.$ This completes the proof of the Theorem.
$\hfill \Box$

\begin{rem} \rm
Note that in proving the above Theorem, we have in fact proved more than advertised.
The next Theorem is actually proved above, but we state it separately for convenience.
\end{rem}

\begin{th}
If $P/Q$ has matrix $$M = M(\vec{a}) = M(a_{1}) \cdots M(a_{n})$$ and  
$R/S$ has matrix $$N = M(\vec{b}) = M(b_{1}) \cdots M(b_{m}),$$
then $\frac{F}{G} = [a_{1},\ldots,a_{n}] \, \sharp[b_{1} \, ,\ldots,b_{m}]$ has matrix $$ M \,  \sharp \, N := M^{T}N',$$ where $N'$ denotes the
matrix obtained by interchanging the rows of $N.$ This gives an explicit formula for
$[a_{1},\ldots,a_{n}] \, \sharp \, [b_{1},\ldots,b_{m}].$ This formula can be used to determine  not only when $N([F/G])$ is unknotted but
also to find its knot type as a rational knot via Schubert's Theorem.
\end{th}

\begin{rem} \rm
Let's draw the exact consequence of Theorem 6. Suppose that 
$$M =
\left( \begin{array}{cc}
P &  Q' \\
Q & U \\
\end{array} \right)$$ and 
$$N =
\left( \begin{array}{cc}
R & S' \\
S & V \\
\end{array} \right).$$
Then
$$\left( \begin{array}{cc}
F & G' \\
G & W \\
\end{array} \right) = M^{T}N' = 
\left( \begin{array}{cc}
P &  Q \\
Q' & U \\
\end{array} \right)
\left( \begin{array}{cc}
S & V \\
R & S' \\
\end{array} \right)$$
$$=
\left( \begin{array}{cc}
PS + QR & PV + QS' \\
Q'S + UR & Q'V + US' \\
\end{array} \right).$$
Thus
$$N([P/Q] + [R/S]) =N([(PS + QR)/(Q'S + UR)])$$ where $|PU - QQ'| = 1.$
\end{rem}

$$ \picill8inby1.8in(UK21) $$
\begin{center}
{\bf Figure 21 - Tangle Concatenation}
\end{center}

 \begin{rem} \rm
In this section we have considered the operation
$$[a_{1},\ldots,a_{n}]\, \sharp \, [b_{1},\ldots,b_{m}] = [a_{n},\ldots,a_{2}, a_{1} + b_{1}, b_{2},\ldots, b_{m}].$$
A simpler operation is given by the following formula
$$[a_{1},\ldots,a_{n}]\,  \&  \, [b_{1},\ldots,b_{m}] = [a_{1},\ldots,a_{n-1}, a_{n}, b_{1}, b_{2},\ldots, b_{m}].$$
This concatenation of continued fractions corresponds directly to the matrix product $M(\vec{a})M(\vec{b}),$ and it has a 
conbinatorial interpretation as a joining of two rational tangles as shown in Figure 21. We remark on this point here to
emphasize that many of the operations in this paper can be systematized in an extended algebra of tangles.  
\end{rem}

 \begin{rem} \rm
Suppose, as in Theorem 5, that $|PS + QR| = 1$ so that
$N([\frac{P}{Q}] + [\frac{R}{S}])$ is an unknot. Then
it follows from Theorem 5 that $N([\frac{Q}{P}] + [\frac{S}{R}])$ is also an unknot. Thus the pair of fractions
$(\frac{P}{Q}, \frac{R}{S})$ yields an unknot if and only if the pair of inverse fractions $(\frac{Q}{P}, \frac{S}{R})$
yields an unknot. In fact, the reader will enjoy proving the following identity for any rational $2$-tangles $A$ and $B:$
$$N(\frac{1}{A} + \frac{1}{B}) = - N(A + B)$$ where the minus sign denotes taking the mirror image of the diagram $N(A + B),$ and equality
is topological equivalence. Thus for any pair of fractions $(\frac{P}{Q}, \frac{R}{S})$ (not necessarily forming an unknot) 
$$N([\frac{P}{Q}] + [\frac{R}{S}]) = - N([\frac{Q}{P}] + [\frac{S}{R}]).$$
\end{rem}

\section{Continued Fractions, Convergents and Lots of Unknots}

Consider a rational fraction, its corresponding continued fraction, and its matrix representation:
$$P/Q = [a_{1},\ldots, a_{n}]$$
with
$$M = M(\vec{a}) = M(a_{1}) \cdots M(a_{n})=
\left( \begin{array}{cc}
P &  Q' \\
Q & U \\
\end{array} \right).$$
Note that since the determinant of this matrix is $(-1)^{n},$ we have the formula
$PU - QQ' = (-1)^{n}$ from which it follows that 
$$P/Q - Q'/U = (-1)^{n}/QU.$$ Hence, by Theorem 5, the diagram $$N([P/Q] - [Q'/U])$$
is unknotted and, as we shall see, is a good candidate to produce a hard unknot.
Furthermore, the fraction $Q'/U$ has an interpretation as the truncation of our
continued fraction $[a_{1},\ldots, a_{n}]:$
$$Q'/U = [a_{1},\ldots,a_{n-1}].$$
To see this formula, let
$$N = M(a_{1}) \cdots M(a_{n-1})=
\left( \begin{array}{cc}
R &  S' \\
S & V \\
\end{array} \right),$$ so that
$$R/S = [a_{1},\ldots, a_{n-1}].$$

Then 
$$\left( \begin{array}{cc}
P &  Q' \\
Q & U \\
\end{array} \right) =
M(a_{1}) \cdots M(a_{n-1})M(a_{n}) = NM(a_{n})$$
$$ = \left( \begin{array}{cc}
R &  S' \\
S & V \\
\end{array} \right)
\left( \begin{array}{cc}
a_{n} &  1 \\
1 & 0 \\
\end{array} \right)$$
$$= \left( \begin{array}{cc}
Ra_{n} + S' &  R \\
Sa_{n} + V & S \\
\end{array} \right).$$
This shows that $Q'/U = R/S = [a_{1},\ldots, a_{n-1}],$ as claimed.
\bigbreak

\begin{defn} \rm
One says that $[a_{1},\ldots, a_{n-1}]$ is a {\it convergent} of $[a_{1},\ldots, a_{n-1}, a_{n}].$
We shall say that two fractions $P/Q$ and $R/S$ are {\it convergents} if the continued fraction of one of them is a convergent of the other.
\end{defn}

\noindent We see from the above calculation that the two consecutive integers $PU$ and $QQ'$ produce two continued fractions
$P/Q = [a_{1},\ldots,a_{n}]$ and $Q'/U = [a_{1},\ldots,a_{n-1}]$ so that the second fraction is a convergent
of the first. 
\bigbreak

\noindent We have proved the following result.
\bigbreak

\begin{th}
Let $P/Q$ and $Q'/U$ be fractions such that the continued fraction of $Q'/U$ is a convergent
of the continued fraction of $P/Q.$ Then $$N([P/Q] - [Q'/U])$$ is an unknot. 
\end{th}

\noindent {\bf Proof.} The proof is given in the discussion above.
\hfill $\Box$ \medskip 

\begin{rem} \rm
This Theorem applies to Figure 5, and our early discussion of the Culprit.
\end{rem}

The property of one fraction being a convergent of the other is in fact, always a property of 
fractions produced from consecutive integers. For example
$$13/9 = [1,2,4] = [1,2,3,1]$$ and 
$$10/7 = [1,2,3]$$ with
$$13 \times 7 - 10 \times 9 = 1,$$
so that $13/9$ and $10/7$ are convergents and 
$N([13/9] + [-10/7])$ is an unknot. We will now prove this statement about convergents. For the next Lemma see also \cite{Ford}.
\bigbreak

\begin{lem}
Let $P$ and $Q$ be relatively prime integers and let $s$ and $r$ be a pair of integers such that
$Ps - Qr = \pm 1.$ Let $R = r + tP$ and $S= s + tQ$ where $t$ is any integer. Then $\{R,S \}$ comprises the set of all
solutions to the equation $PS - QR = \pm 1.$ If $Ps - Qr = \pm 1$ and $PS - QR = \mp 1,$ Then all solutions are given in the form
$R = -r + tP$ and $S= -s + tQ.$
\end{lem}

\noindent {\bf Proof.} Without loss of generality we can assume that $Ps - Qr = 1.$ We leave it to the reader to 
formulate the case where $Ps - Qr= -1.$
Certainly $R$ and $S$ as given in the statement of the Lemma satisfy the equation $PS - QR = 1.$
Suppose that $R$ and $S$ is some solution to this equation. Then it follows by taking the difference with the equation
$Ps - Qr = 1$ that $P(S-s) = Q(r - R).$ Since $P$ and $Q$ are relatively prime it follows at once from this equation and 
the uniqueness of prime factorization of integers that  $R = r + tP$ and $S=s + tQ$ for some integer $t.$ The last part of the Lemma follows 
by the same form of reasoning. This completes
the proof of the Lemma.
\hfill $\Box$ \medskip 

\begin{th}
Let $P$ and $Q$ be relatively prime integers and let $P/Q = [a_{1},\ldots,a_{n}]$ be a 
continued fraction expansion for $P/Q.$ Let $r/s = [a_{1},\ldots,a_{n-1}]$ be the convergent for $[a_{1},\ldots,a_{n}].$
Let $R = r + tP$ and $S= s + tQ$ where $t$ is any integer. Then $R/S = [a_{1},\ldots,a_{n}, t].$ Thus $P/Q$ is a convergent
of $R/S.$ We conclude that if $P/Q$ and $R/S$ satisfy the condition that $N([P/Q] - [R/S])$ is an unknot, then one
of $P/Q$ and $R/S$ is a convergent of the other.
\end{th}

\noindent {\bf Proof.} By the previous discussion and the hypothesis of the Theorem, we are given that 
$$M(a_{1}) \cdots M(a_{n}) = 
\left( \begin{array}{cc}
P & r \\
Q & s \\
\end{array} \right),$$
and
$$\left( \begin{array}{cc}
P &  r \\
Q & s \\
\end{array} \right)
\left( \begin{array}{cc}
t & 1 \\
1 & 0 \\
\end{array} \right) =
\left( \begin{array}{cc}
Pt + r &  P \\
Qt + s & Q \\
\end{array} \right).$$
By the previous Lemma the set $\{ R=r + tP, \, S=s + tQ \}$ comprises all solutions to 
$PS - QR = (-1)^{n}.$ Thus all such solutions correspond to continued fractions of which $P/Q$ is a convegent.
This completes the proof of the Theorem.
\hfill $\Box$ \medskip

\noindent{\bf Infinite Continued Fractions.}
Note that if we start with an infinite continued fraction
$$[a_{1}, a_{2}, a_{3}, \ldots ]$$ and define
$$P_{n}/Q_{n} = [a_{1},\ldots, a_{n}],$$ then we have
$$P_{n}/Q_{n} - P_{n-1}/Q_{n-1} = (-1)^{n}/Q_{n}Q_{n-1}.$$
This equation is usually used to analyze the convergence of the trucations of the continued fraction to a real number.
Here we see the single infinite continued fraction, whether or not convergent, producing an infinite family of unknot
diagrams: $$K_{n} = N([P_{n}/Q_{n}] - [P_{n-1}/Q_{n-1}]).$$
\bigbreak

Surely, the simplest instance of this phenomenon occurs with the infinite continued fraction for the Golden Ratio:
$$\frac{1+\sqrt{5}}{2} = [1,1,1,1,1,1,\ldots].$$
Here the convergents are $f_{n}/f_{n-1} = [1,1,1,\ldots,1]$ (with $n$ $1$'s) and $f_{n}$ is the $n$-th Fibonacci
number where $f_{-1} = 0, f_{0} = 1, f_{1} = 1$ and $f_{n+1} = f_{n} + f_{n-1}$ for $n \ge 1.$ The matrices are
$$M(1)^{n} =
\left( \begin{array}{cc}
f_{n} &  f_{n-1} \\
f_{n-1} & f_{n-2}\\
\end{array} \right).$$ Thus we have
$$f_{n}/f_{n-1} - f_{n-1}/f_{n-2} = (-1)^{n}/f_{n-1}f_{n-2}.$$ In this instance, this formula shows that the truncations 
converge to the Golden Ratio.
\bigbreak

In Figure 22 we illustrate a hard unknot based on the eight and seventh truncations of the Golden Ratio continued fraction.
This example is produced by using the standard closure of 
$N([f_{8}/f_{7}] - [f_{7}/f_{6}]( = N([34/21] - [21/13])$. In Figure 23 we illustrate the closure of 
$N([34/21] - [21/13])$ once again, but this time we use the first tangle turned by $180$ degrees. The result is a hard knot
diagram when we tuck an arc across the closure. In Section 8 we discuss the general constructions behind these two cases.

$$ \picill3.7inby2in(UK22) 
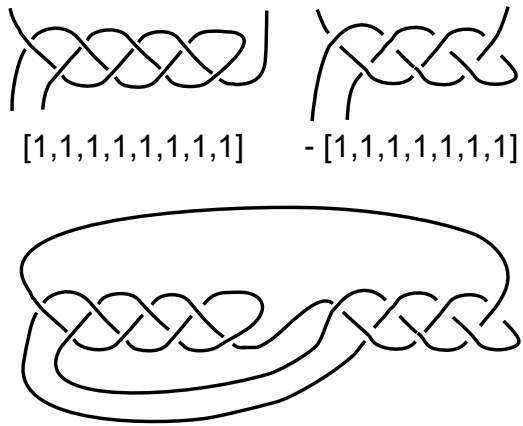
$$
\begin{center}
{\bf Figure 22 - $N([f_{8}/f_{7}] - [f_{7}/f_{6}]) = N([34/21] - [21/13])$}
\end{center}

$$ \picill3.8inby4.3in(UK23) $$
\begin{center}
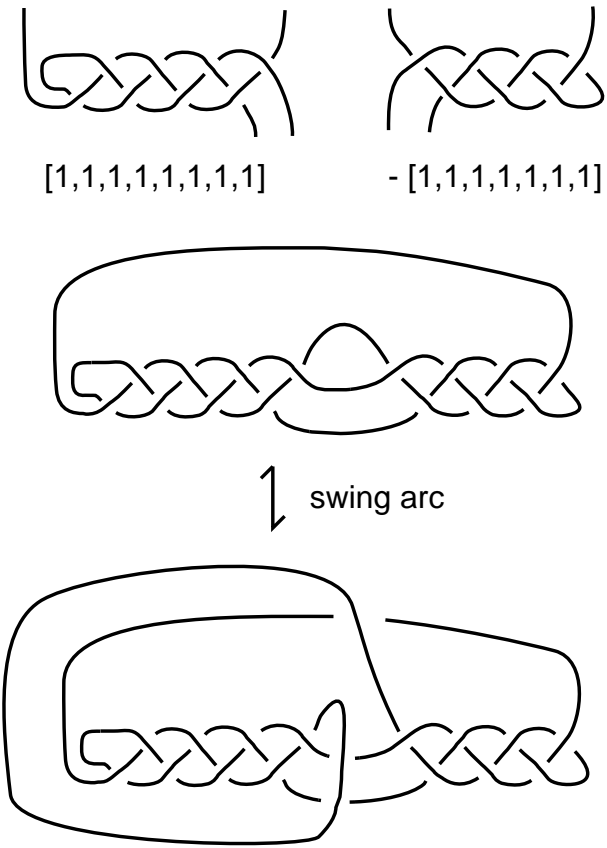
{\bf Figure 23 - $N([f_{8}/f_{7}] - [f_{7}/f_{6}])$ reversed and tucked}
\end{center}

\section{Constructing Hard Unknots}

In this section we indicate how to construct hard unknots by using positive alternating tangles $A$ and $B$ such that
$N(A - B)$ is unknotted. By our main results we know how to construct infinitely many such pairs of tangles by taking 
a continued fraction and its convergent, with the corresponding tangles in reduced (alternating) form.

$$ \picill5inby2.1in(UK24) $$
\begin{center}
{\bf Figure 24 - $K = N([1,4] - [1,3])$ }
\end{center}

Let's begin with the case of $5/4 = [1,4] = [1,3,1]$ and $4/3 = [1,3].$ In Figure 24 we show the standard representations of 
$[1,4]$ and $[1,3]$ as tangles, and the corresponding construction for the diagram of $K = N([1,4] - [1,3]).$ The reader will
note that this diagram is a hard unknot with 9 crossings, one less than our original Culprit of Figure 5. We give another version
of it in Figure 31 (equivalent to its mirror image $H$). In Section 6 we show that $H$
is one of a small collection of minimal hard unknot diagrams having the form $N(A - B)$ for reduced positive rational tangle diagrams $A$ and
$B$.

$$ \picill5inby2in(UK25) $$
\begin{center}
{\bf Figure 25 - $N([1,3] - [1,2])$}
\end{center}

In most cases, if one takes the standard representations of the tangles $A$ and $B$, and forms the diagram for
$N(A - B)$, the resulting unknot diagram will be hard. There are some exceptions however, and the next example illustrates 
this phemomenon. 
\smallbreak

\noindent In Figure 25 we show the standard representations of 
$[1,3]$ and $[1,2]$ as tangles, and the corresponding construction for the diagram of $N([1,3] - [1,2]).$ This diagram, while
unknotted, is not a hard unknot diagram due to the three-sided outer region. This outer region allows a type III Reidemeister
move on the surface of the two dimensional sphere. In this example, tucking an arc does not create a hard unknot from the given diagram
(there is be a type III move available after the tuck).

$$ \picill5inby3.5in(UK26) $$
\begin{center}
{\bf  Figure 26 -$N([1,3] - [1,2]^{vflip})$ }
\end{center}

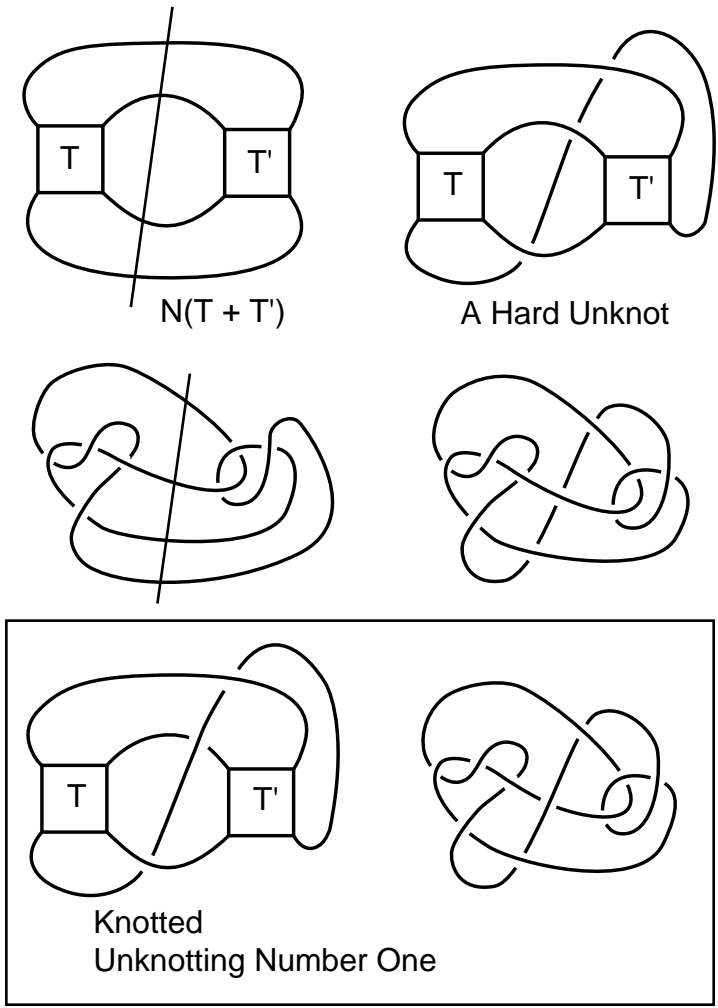
\noindent {\bf The Tucking Construct.} Figure 26 shows a way to remedy this situation. Here we have replaced $[1,2]$ by
$[1,2]^{vflip},$ the
$180$ degree turn of the  tangle $[1,2]$ about the vertical direction in the page. Now we see that the literal diagram of $N([1,3] -
[1,2]^{vflip})$ is of course still unknotted and is also  not a hard unknot diagram. However this diagram can be converted to an unknot diagram
by tucking an arc as shown in the  Figure. The resulting hard unknot is the same diagram of $10$ crossings that we had in Figure 5 as our
initial Culprit. Note that the other possibility of flipping both tangles in Figure 25 or flipping the first tangle do not lead to
a hard unknots. We call this strategem the {\it tucking construct}. Tucking is accompanied by the vertical flip of on one of the
tangles to avoid  the placement of a Reidemeister move of type III as a result of the tuck.
\bigbreak

\noindent {\bf The Culprit Revisited}
Let's consider the example in Figure 5 again. Here we have $P/Q =F(A) = -3/4$ and $R/S = F(B) = 2/3.$
We have $P/Q + R/S = -3/4 + 2/3 = -1/12.$ Thus $N([-3/4] + [2/3])$ is an unknot by Theorem 5. This is exactly the unknot
$C'$ illustrated in  Figure 5.

$$ \picill4inby5.4in(UK27) $$
\begin{center}
{\bf Figure 27 - The tucking construct}
\end{center}

We can make infinitely many examples of this type. View Figure 27. The pattern is as follows.
Suppose that $T = [P/Q]$ and $T' = [R/S]$ are rational tangles such that $PS - QR = \pm 1.$ Then we know that
$N(T - (T')^{vflip})$ is an unknot. Furthermore we can assume that each of the tangles $T$ and $T'$ are in alternating form.
The two tangle fractions have opposite sign and hence the alternation of the weaves in
each tangle will be of opposite type. We create a new diagram for $N(T - (T')^{vflip})$ by putting an arc from the bottom of the closure 
entirely underneath the diagram as shown in Figure 27. This is an example of a sucessful tucking construct. Note how in the
example shown in Figure 27, the knot diagram resulting from the tucking construction is indeed our original hard unknot diagram.
There are no  simplifying Reidemeister moves and there are no moves of type III available on the diagram.

$$ \picill1.8inby1.9in(UK28) $$
\begin{center}
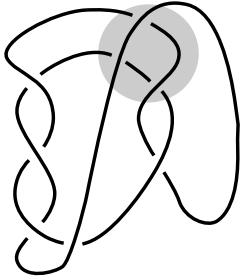
{\bf Figure 28 - The tucking construct can produce an easy unknot}
\end{center}

$$ \picill3.3inby1.7in(UK29) $$
\begin{center}
{\bf Figure 29 -Another hard unknot}
\end{center}

This process will always produce 
hard unknot diagrams, except in some easy special cases. View, for example, Figure 28. We leave investigation of these cases as an
exercise for the  reader.
\bigbreak

View Figure 29. Here we illustrate a typical hard unknot produced by the tucking construction.
In this case, once the arc is swung out, one sees that the diagram is the numerator of the sum of the tangles
$$8/3 = [2,1,2]$$ and 
$$-13/5 = -[2,1,1,1,1]$$ with
$$8/3 - 13/5 = 1/15.$$

\begin{rem} \rm
View the bottom of Figure 27 and note that we have indicated that the hard unknot becomes an alternating
and hence knotted (this can be proved by applying the results in \cite{KState}) diagram by switching one crossing in the tucking construction.
Thus the tucking construction is also an infinite source of knots of unknotting number one.
\end{rem}

$$ \picill2.5inby1.6in(UK30) $$
\begin{center}
{\bf Figure 30 - A complex unknot from NotKnot}
\end{center}

It is interesting to comb the literature for unknots and see that many of them are produced by the tucking construct.
For example, in Figure 30 we illustrate a knot from page 5 of the ``Supplement to Not Knot" by David Epstein and Charlie Gunn
\cite{NotKnot}.
In this case we see that their example can be seen as the tucking construct applied to the rational tangles $[11/7]$ and
$[-8/5].$ The diagram itself is not a hard unknot, since there are Reidemeister III moves available. However, it is not difficult to see that,
using Reidemeister moves, this diagram must be made more complex before it simplifies to the unknot. This example shows that we could
consider a wider class of unknotted diagrams than those that have no Reidemeister III moves and no simplifying moves. We can call a diagram
in this wider a class a {\it complex} unknot diagram if it requires diagrams with more crossings in order to be unknotted. The KnotNot diagram
of Figure 30 is complex.

\section{The Smallest Hard Unknots}

Figure 31 illustrates two hard unknot diagrams $H$ and $J$ with $9$ crossings.
\bigbreak

\begin{conj}
Up to mirror images and flyping tangles in the diagrams, the hard unknot diagrams $H$ and $J$ of 9 crossings, shown in Figure 31
($K = -H$ appears earlier in Figure 24), have the least number of crossings among all hard unknot diagrams.
\end{conj}

$$ \picill2inby3.3in(UK31) $$
\begin{center}
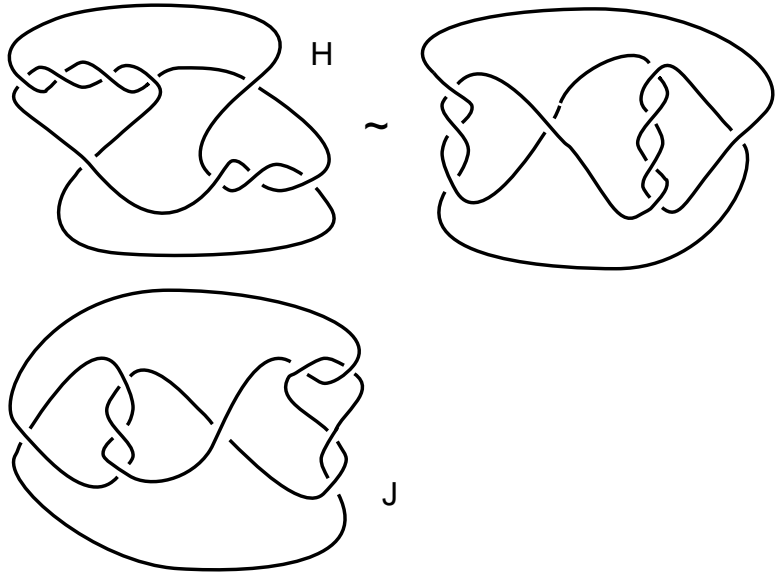
{\bf Figure 31 - $H$ and $J$ are hard unknots of $9$ crossings}
\end{center}

Two equivalent versions of the diagram $H$ appear in Figure 31. The right-hand version of $H$ in this figure
is of the form 
$$H = N([1 + 1/3] - [1 + 1/4]) = N([1,3] -[1,4]) = N([4/3] - [5/4]).$$ 
Note that $[1,3]$ and $[1,4]=[1,3,1]$ are convergents.
Note also that the diagram $K$ of Figure 24 is given by 
$K = N([1,4] - [1,3]) = -N([1,3] - [1,4]) = -H.$ Thus $H$ and $K$ are mirror images of each other.
\bigbreak

\noindent The diagram $J$ in Figure 31 is of the form $$N([1 + 1/3] - [1 + 1/(2 + 1/2)]) = N([1,3] - [1,2,2]) = N([4/3] -
[7/5]).$$ Note that $[1,3]=[1,2,1]$ and $[1,2,2]=[1,2,1,1]$ are convergents.
\bigbreak

\noindent Note also that the crossings in $J$ corresponding to $1$ in $[1,3]$ and $-1$ in $-[1,2,2]=[-1,-2,-2]$ can be switched and we will
obtain another diagram $J',$ arising as sum of two alternating rational tangles, that is also a hard unknot. This diagram can be obtained from
the diagram $J$ without switching crossings by performing flypes (Figure 13) on the subtangles $[1,3]$ and $[1,2,2]$ of $J,$ and then doing an
isotopy of this new diagram on the two dimensional sphere. (We leave the verification of this statement to the reader.) Thus the diagram
$J'$ can be obtained from $J$ by flyping. 
A similar remark applies to the diagram $H,$ giving a corresponding diagram $H',$ but in this case $H'$ is easily seen to be equivalent to $H$
by an isotopy that does not involve any Reidemeister moves. Thus, up to these sorts of modifications, we have produced essentially two hard
diagrams with $9$ crossings. Other related hard unknot diagrams of $9$ crossings can be obtained from these by taking mirror images.
\bigbreak

\noindent We have the following result.

\begin{th}
The diagrams $H$ and $J$  shown in Figure 31 are, up to flyping subtangle diagrams and taking mirror images, the smallest
hard unknot diagrams in the form 
$N(A - B)$ where $A$ and $B$ are rational tangles in reduced positive alternating form.
\end{th}

\noindent {\bf Proof.} It is easy to see that we can assume that $A=[P/Q]$ where $P$ and $Q$ are positive, relatively 
prime and $P$ is greater than $Q.$ We leave the proof that one can choose $P$ greater than $Q$ to the reader, with the hint:
Verify that the closure diagram in Figure 25 is equivalent to the diagram in Figure 31 on the surface of the two dimensional 
sphere, without using any Reidemeister moves. 
\bigbreak

We then know from Theorem 8 that $B=[-R/S]$ where one of $P/Q$ and $R/S$  is a convergent of 
the other. We can now enumerate small continued fractions. We know the total of all terms in $A$ and $B$ must be less than or
equal to $9$ since $H$ and $J$ each have nine crossings.
\bigbreak

In order to make a $9$ crossing unknot example of the form $N(A - B)$ where $A$ and $B$ are rational tangles in reduced positive alternating
form, we must partition the number $9$ into two parts corresponding to the number of crossings in each tangle. It is not hard to see that
we need to use the partition $9 = 4 + 5$ in order to make a hard unknot of this form. Furthermore, $4$ must correspond to the the 
continued fraction $[1,3],$ as $[2,2]$ will not produce a hard unknot when combined with another tangle. Thus, for producing $9$ crossing
examples we must take $A=[1,3].$ Then, in order that $A$ and $B$ be convergents, and $B$ have $5$ crossings, the only possibilities for
$B$ are $B=[1,4]$ and $B=[1,2,2].$ These choices produce the diagrams $H,H',J,J'.$ It is easy to see that no diagrams with less than 
$9$ crossings will suffice to produce hard unknots, due to the appearance of Reidemeister moves related to the smaller partitions.
This completes the proof.
$\hfill \Box$

\section{The Goeritz Unknot}

The earliest appearance of a hard unknot is a 1934 paper of Goeritz \cite{Goeritz}. In this paper Goeritz gives the hard
unknot shown in Figure 32. As the reader can see (for example by twisting vertically the tangle [-3] twice), this example is a
variant on $N([4] + [-3])$ which is certainly unknotted. The Goeritz example $G$ has 11 crossings, due to the extra two twists
that make it a hard unknot. It is part of an infinite family based on $N([n] +[-n+1]).$

$$ \picill2inby1.2in(UK32) $$
\begin{center}
{\bf Figure 32 - The Goeritz Hard Unknot}
\end{center}

$$ \picill8inby4.1in(UK33) $$
\begin{center}
{\bf Figure 33 - Unknotting the Goeritz Unknot}
\end{center}

In Figure 33 we illustrate steps in an unknotting process for $G.$
The moves labeled {\em swing} consist in swinging the corresponding marked arc in the diagram.
These moves when realized by Reidemeister moves will complicate the diagram before taking it to the indicated simplification.
Notice that after two swing moves the Goeritz diagram has been transformed into the diagram $H$ that is again a hard unknot 
diagram, but with 9 crossings rather than 11! This is the small hard unknot diagram that we discussed in the
previous section. Figure 33 then  continues, with one more swing move applied to $H$ (as is necessary) and
then a sequence of Reidemeister III moves, a Reidemeister I move and a II move, taking the diagram to a simple twist that
unknots by a few more simplifying  Reidemeister moves. This example is instructive in that it shows that a sequence of swing moves
may be needed to undo a hard unknot and that the number of crossings may go up and down in the process.

\section{Stability in Processive DNA Recombination}

In this section we use the techniques of this paper to study properties of processive DNA recombination topology.
Here we use the tangle model of DNA recombination \cite{Sumners1,Sumners2,Su} developed by C. Ernst and D.W. Sumners. In this model the DNA is
divided into two regions corresponding to  two tangles $O$ and $I$ and a recombination site that is associated with $I.$ This division is a
model of how the  enzyme that performs the recombination traps a part of the DNA, thereby effectively dividing it into the tangles $O$ and $I.$ 
The recombination site is
represented by another tangle $R.$ The entire arrangement is then a knot or link $K[R] = N(O + I + R).$ We then consider a single
recombination in the form of starting with $R = [0],$ the zero tangle, and replacing $R$ with the tangle $[1]$ or the tangle
$[-1].$ Processive recombination consists in consecutively replacing again and again by $[1]$ or by $[-1]$ at the same site. Thus, in processive
recombination we obtain the knots and links $$K[n] = N(O + I + [n]).$$ The knot or link $K[0] = N(O + I)$  is called the  {\em DNA substrate},
and the tangle $O + I$ is called the {\em substrate tangle}. It is of interest to obtain a uniform formula for knots and links $K[n]$ that
result from the processive recombination.
\bigbreak

In some cases the substrate tangle is quite simple and is represented as a single tangle $S = O + I.$ For example, we illustrate
processive recombination in Figure 34 with $S = [-1/3] = [0,-3]$ and $I = [0]$ with $n = 0,1,2,3,4.$ Note that by Proposition 1 of
Section 2.1,
$$K[n] = N(S + [n]) = N([0,-3] + [n]) = N([-3, 0 + n]) = N([-3,n])$$
$$= N(-[3,-n]) = N(-[2,1,n-1]).$$ This formula gives the
abstract form of all the knots and links that arise from this recombination process. We say that the formula $$K[n] =
N(-[2,1,n-1])$$ for
$n>1$ is {\it stabilized} in the sense that all the terms in the continued fraction have the same sign and the $n$ is in one
single place in the fraction. In general, a {\it stabilized fraction} will have the form
$$N(\pm [a_{1},a_{2},\ldots a_{k-1}, a_{k} + n, a_{k+1}, \ldots ,a_{n}])$$
where all the terms $a_i$ are positive for $i\ne k$ and $a_k$ is non-negative.

$$ \picill5inby6in(UK34) $$
\begin{center}
{\bf Figure 34 - Processive Recombination with $S = [-1/3]$}
\end{center}

Let's see what the form of the processive recombination is for an arbitrary sequence of recombinations. We start with 
$$O = [a_1, a_2 , \ldots, a_{r-1}, a_{r}]$$
$$I = [b_1, b_2, \ldots, b_{s-1}, b_{s}].$$
Then
$$K[n] = N(O + (I + [n])) = N([a_1, a_2 , \ldots, a_{r-1}, a_{r}] + [n + b_1, b_2, \ldots, b_{s-1}, b_{s}])$$
$$K[n] = N([a_{r}, a_{r-1} , \ldots, a_{2}, a_1 + n + b_1, b_2, \ldots, b_{s-1}, b_{s}]).$$
\bigbreak

\begin{prop}
The formula
$$K[n] = N([a_{r}, a_{r-1} , \ldots, a_{2}, a_1 + n + b_1, b_2, \ldots, b_{s-1}, b_{s}])$$
can be simplified to yield a stable formula for the processive recombination when $n$ is
sufficiently large. 
\end{prop}

\noindent {\bf Proof.}  Apply Proposition 1 of Section 2.1.
$\hfill \Box$
\bigbreak

Here is an example. Suppose we take $O = [1,1,1,1]$ and $I=[-1,-1,-1]$ so that the DNA substrate is an 
(Fibonacci) unknot. ($I$ is the negative of the convergent of $O.$) Then by the above calculation
$$K[n] = N([1,1,1,1 + n + (-1),-1,-1]) = N([1,1,1,n,-1,-1]).$$
Suppose that $n$ is positive.  
Applying the reduction formula of Proposition 1, we get
$$K[n] = N([1,1,1,n,-1,-1]) = N([1,1,1,n-1,1,0,1]) = N([1,1,1,n-1,2]),$$ and this is a stabilized form for the processive 
recombination.
\bigbreak

More generally, suppose that $O=[a_{1},a_{2},\ldots,a_{n}]$ where all of the $a_{i}$ are positive. Let
$I = [-a_{1},-a_{2},\ldots,-a_{n-1}].$ Then $K[0] = N(O + I)$ is an unknotted substrate by our result about
convergents. Consider $K[n]$ for positive $n.$ We have
$$K[n] = N([a_{n},a_{n-1},\ldots,a_{2}, a_{1} + n - a_{1}, -a_{2},\ldots,-a_{n-1}])$$
$$= N([a_{n},a_{n-1},\ldots,a_{2}, n, -a_{2},\ldots,-a_{n-1}])$$
$$= N([a_{n},a_{n-1},\ldots,a_{2}, (n-1), 1,  (a_{2}-1), a_{3},\ldots, a_{n-1}]).$$
If $a_{2}-1$ is not zero, the process terminates immediately. Otherwise there is one more step.
In this way the knots and links proceeding from the recombination process all have a uniform stabilized form.
Further successive recombination just adds more twist in one entry in the continued fraction diagram whose closure is
$K[n].$
\bigbreak

The reader may be interested in watching a visual demonstration of these properties of DNA recombination.
For this, we recommend the program {\it Ginterface} (TangleSolver) \cite{Vasquez} of Mariel Vasquez. Her program
can be downloaded from the internet as a Java applet, and it performs and displays DNA recombination.
Figure 35 illustrates the form of display for this program. The reader should be warned that the program uses
the reverse order from our convention when listing the terms in a continued fraction. Thus we say $[1,2,3,4]$ 
while the program uses $[4,3,2,1]$ for the same structure.

$$ \picill5inby3in(UK35.EPSF) $$
\begin{center}
{\bf Figure 35 - Processive Recombination with $S = [1,1,1,1] + [-1,-1,-1]$}
\end{center}

\section{Collapsing Tangles and Hard Knots}

In this section we point out that our results about when the closure of a sum of rational tangles is unknotted can be 
generalized to collapses of tangle sums to particular knots and links.
\bigbreak

\noindent {\bf Notation.} Let $R$ be any tangle.  We shall  write $$S=[a_1, \ldots ,a_n, R]$$
to denote the tangle $$[a_1] + \frac{1}{[a_2] +\frac{1}{ \ldots +\frac{1}{[a_n] + \frac{1}{R}}}}.$$
That is, we will occasionally replace $[a]$ by $a$ in the notation for continued fraction forms of tangles. 
\bigbreak

View Figure 36. This Figure illustrates the identity
$$N(\frac{1}{[m] + \frac{1}{[n] + T}} - \frac{1}{[m]}) = D(T)$$
where $D(T)$ denotes the denominator of $T,$ and $T$ is any tangle.
The example is due to DeWitt Sumners and is fashioned to illustrate that an equation
$$N(A + B) = K$$ for fixed $K$ and variable $A$ and $B$ may have infinitely many solutions.
So far in this paper we have concentrated on the case of this phenomenon where $K$ is an unknot.
In point of fact, for any knot $K$ this equation has infinitely many solutions beyond the examples just 
indicated. But first, let's analyse the class of examples in Figure 36.

$$ \picill8inby4.3in(UK36) $$
\begin{center}
{\bf Figure 36 - Sumners' Example}
\end{center}

\noindent Here is the analysis of Figure 36 in our formalism.
Note that for any tangle $T$, 
$$N([T,0]) = N(T + 1/[0]) =  N(T + [\infty]) = D(T).$$
We have

$$N(\frac{1}{[m] + \frac{1}{[n] + T}} - \frac{1}{[m]}) = N([0,m,[n]+T]+[0,-m]) = $$
$$N([[n]+T, m, 0 + 0, -m]) = N([n+T, m, 0 , -m]) =$$
$$N([[n] + T, 0]) = D([n] + T) = D(T).$$  
\bigbreak

\noindent More generally we have the following Proposition, whose proof we leave to the reader.
\smallbreak

\begin{prop}
Let $T$ and $S$ denote any tangles. Then 
\begin{enumerate}
\item $N([a_{1},a_{2},\ldots ,a_{n-1}, a_{n} + T] - [a_{1},a_{2},\ldots ,a_{n-1}]) = D(T),$ 
\item $N([a_{1},a_{2},\ldots ,a_{n-1}, a_{n} + T] - [a_{1},a_{2},\ldots ,a_{n}]) = N(T),$
\item $N([a_{1},a_{2},\ldots ,a_{n-1}, a_{n} + T] - [a_{1},a_{2},\ldots ,a_{n-1},a_{n} - S]) = N(T + S).$
\end{enumerate}
\end{prop}

\begin{rem} \rm
Note that $N(T + [\infty]) = D(T)$ for any tangle $T.$
Note that $N([\infty])$ is the unknot, and that since $1/[\infty] = [0],$ we have the tangle equations
$$[a_{1},a_{2},\ldots ,a_{n-1}, a_{n} + [\infty]] = [a_{1},a_{2},\ldots ,a_{n-1}, [\infty]] = [a_{1},a_{2},\ldots ,a_{n-1}].$$
Hence if $T = [\infty],$ then the second part of this Proposition shows that 
$$N([a_{1},a_{2},\ldots ,a_{n-1}] - [a_{1},a_{2},\ldots ,a_{n}])$$ is unknotted, as we have shown in Theorem 8.
In the third case we see that the numerator of the sum of any two tangles can be separated by rational twisting that will cancel
back to the original sum. All rational knots are susceptible to such decompositions by using rational tangles $T$ and $S.$
This Proposition gives the basic generalization of our unknotting result to general tangle collapse to a knot or link of the
form $N(T + S)$. For example, we have that $N([1,1,1,2] + [-1,-1,-1,2]) = N([2,2]).$ This is a collapse to the figure eight
knot. We trust that the reader will enjoy making other examples to illustrate this result.
\end{rem}

\section{Recalcitrance Revisited}

$$ \picill3inby2.2in(UK37) 
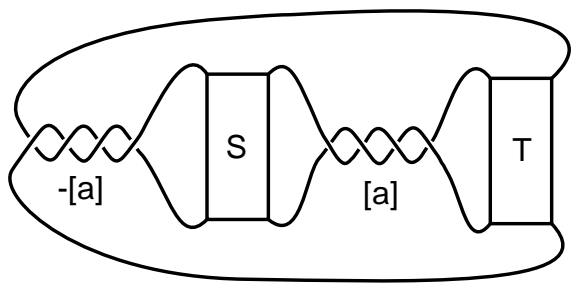
$$
\begin{center}
{\bf Figure 37 - Twisting Up the Recalcitrance}
\end{center}

Suppose that $K[a]$ denotes the knot diagram shown in Figure 37, so that 
$$K[a] =N(-[a] + [S] + [a] + T]),$$ and we shall assume that $K[0] = N(S + T)$ is unknotted.
We assume that the diagram $K[a]$ is hard, a generalization of the form of the hard diagram $H$ in Figure 25.
Let's also assume that it takes $N$ Reidemeister moves to transform 
$[-1] + S + [1]$ to $S.$ This transformation is easily accomplished in three-dimensions by one full rotation, but
may require many Reidemeister moves in the plane (keeping the ends of the tangle fixed). It is certainly possible to produce
alternating rational tangles $S$ and $T$ with this property. Let $C = C(S) + C(T)$ denote the total number of crossings in the tangles
$S$ and $T.$ Then the recalcitrance (see Section 1) of $K[0]$ is $R(0) = k/C$ for some $k$ and the recalcitrance of 
$K[a]$ is generically given by the formula $$R(a) = \frac{aN + k}{2a + C}$$ since each extra turn of $S$ will add $N$
Reidemeister moves to the untying, and the number of crossings of $K[a]$ is equal to $2a + C.$ We conclude from this
that {\it for large $a$ the recalcitrance $R(a)$ is as close as we like to $N/2.$} This shows that {\it the recalcitrance of an
unknotted diagram can be arbitrarily large.} There is no upper limit for the ratio of the number of moves needed to undo the knot
in relation to the number of crossings in the original knot diagram.

\section{Afterthoughts - Farey Series, Continued Fractions, Pick's Theorem and Ford Circles}

The theme of this paper has been our result (Theorem 5) that given two fractions $\frac{p}{q}$ and $\frac{r}{s}$ such that 
$|ps - qr| = 1,$ we can construct an unknot diagram $N([\frac{p}{q}] - [\frac{r}{s}])$ from the rational tangles $[\frac{p}{q}]$ and
$[\frac{r}{s}]$. The arithmetic  condition \,\,\,\, $|ps - qr| = 1$ has so many beautiful and surprising connections to other mathematics, that
we must mention them in this last section of the paper! We shall touch on Farey series \cite{Farey}, the Riemann Hypothesis, continued
fractions, Pick's Theorem  and  Ford circles \cite{C2}. There is more, but we hope that this will give the reader a
taste, and perhaps there will arise  new connections with knots and unknots as well through this discussion. 
\bigbreak

\noindent To remind us to  think of knots, we shall call fractions $\frac{p}{q}$ and $\frac{r}{s}$ such that 
$|ps - qr| = 1$ an {\it unknot pair}. Note that if $\frac{p}{q}, \frac{r}{s}$ is an unknot pair, then $\frac{q}{p}, \frac{s}{r}$ 
is also an unknot pair (see Remark 4 in Section 3). When $\frac{p}{q} < \frac{r}{s},$ then we have $ps - qr = -1,$ the plus sign
appearing when the  fractions are in the other order. It is convenient to include the formal fractions $0/1$ and $1/0$ in these discussions
just as we have  done earlier in the paper with the tangles $[0]$ and $[\infty].$

\subsection{Farey Series and Continued Fractions}

Given two fractions $\frac{a}{b}$ and $\frac{c}{d}$ such that $ad - bc = -1$, one can form the {\it mediant} fraction
$\frac{e}{f} = \frac{a+c}{b+d}.$  If $\frac{a}{b} < \frac{c}{d}$, then it is easy to see that 
$\frac{a}{b} < \frac{e}{f} < \frac{c}{d}$ and that $af - be = -1 = ed - fc.$ This means that one can iterate the mediant construction,
producing infinitely many unknot pairs from the given pair $(\frac{a}{b}, \frac{c}{d})$. 
\bigbreak

\noindent This iteration via the mediant fraction was discovered by John Farey \cite{Farey} in 1816. It follows from Farey's construction
that if one starts with the fractions  $(\frac{0}{1}, \frac{1}{1})$ and iterates the mediant construction forever, then all rational
numbers in the interval $[0,1]$ will appear uniquely in their reduced forms. 
\bigbreak

\noindent In order to create any positive
rational number greater than $1$ we also consider the formal fraction $\frac{1}{0}$ and
apply the mediant construction. Note that, then $\frac{1}{1} =\frac{0+1}{1+0}$ (one is the mediant of zero and infinity) so $\frac{0}{1}$ and
$\frac{1}{0}$ can be seen as the ``primal ancestors" of the rational numbers.
\bigbreak

\noindent Given a real number $x>1$, one can consider the rational numbers which, when expressed in lowest terms, have denominators
less than $x.$ The {\it Farey series} corresponding to $x$ is the (ordered) set of positive rational numbers less than  or equal to $1$ which,
in reduced form, have denominators less than $x.$  For example, the Farey series corresponding to $6$ is (in order):
$$\{\frac{1}{5},\,\,\,\frac{1}{4},\,\,\,\frac{1}{3},\,\,\,\frac{2}{5},\,\,\,\frac{1}{2},\,\,\,\frac{3}{5},\,\,\,\frac{2}{3},\,\,\,\frac{3}{4},\,\,\,\frac{4}{5},\,\,\,\frac{1}{1}\}.$$
The reader will notice that each adjacent pair of fractions in the Farey series above is an unknot pair.
Each Farey series can be obtained by
iterating the mediant construction starting with $\{\frac{0}{1}, \frac{1}{1}\},$ and omitting at each step the
fractions whose denominators exceed $x$, until we reach a set for which every
mediant is no longer permitted.
For example, here is the generation of the Farey series
corresponding to $6$ (with $\frac{0}{1}$ and $\frac{1}{1}$ retained at the left and the right).
\bigbreak

$$\{\frac{0}{1},\,\,\,\, \frac{1}{1}\}$$
$$\{\frac{0}{1},\,\,\frac{1}{2},\,\, \frac{1}{1}\}$$
$$\{\frac{0}{1},\,\,\frac{1}{3},\,\,\,\frac{1}{2},\,\,\,\frac{2}{3},\,\,\frac{1}{1}\}$$
$$\{\frac{0}{1},\,\,\frac{1}{4},\,\,\,\frac{1}{3},\,\,\,\frac{2}{5},\,\,\,\frac{1}{2},\,\,\,\frac{3}{5},\,\,\,\frac{2}{3},\,\,\,\frac{3}{4},\,\,\,\,\,\frac{1}{1}\}$$
$$\{\frac{0}{1},\,\,\frac{1}{5},\,\,\,\frac{1}{4},\,\,\,\frac{1}{3},\,\,\,\frac{2}{5},\,\,\,\frac{1}{2},\,\,\,\frac{3}{5},\,\,\,\frac{2}{3},\,\,\,\frac{3}{4},\,\,\,\frac{4}{5},\,\,\frac{1}{1}\}$$
$$\{\frac{0}{1},\,\,\frac{1}{6},\,\,\frac{1}{5},\,\,\,\frac{2}{9},\,\,\frac{1}{4},\,\,\,\frac{2}{7},\,\,\frac{1}{3},\,\,\,\frac{3}{8},\,\,\frac{2}{5},\,\,\,\frac{3}{7},\,\,\frac{1}{2},\,\,\,\frac{4}{7},\,\,\frac{3}{5},\,\,\,\frac{5}{8},\,\,\frac{2}{3},\,\,\,\frac{5}{7},\,\,\frac{3}{4},\,\,\,\frac{7}{9},\,\,\frac{4}{5},\,\,\frac{5}{6},\,\,\frac{1}{1}\}$$
Adjacent pairs are unknot pairs because each fraction is the mediant of its neighbors.
\bigbreak

In fact, given a reduced positive fraction of the form $p/q$ with $0 < p < q$ we can find the two neighboring fractions that give rise to it in the Farey process above. We will give the method for doing this now. Let 
$$p/q = [0,a_{1},a_{2},\cdots, a_{n}]$$ be the positive continued fraction expression for
 $p/q.$ Note that $[0,a_{1},a_{2},\cdots, a_{n}]$ is uniquely determined except for the last term.
 If $a_{n}=1,$ then we can rewrite $$[0,a_{1},a_{2},\cdots, a_{n-1},1]=[0,a_{1},a_{2},\cdots, a_{n-1}+1],$$
 and if $a_{n} >1,$ then we can rewrite 
  $$[0,a_{1},a_{2},\cdots, a_{n}]=[0,a_{1},a_{2},\cdots, a_{n}-1,1].$$ Thus it suffices to assume that we deal with $$p/q = [0,a_{1},a_{2},\cdots, a_{n}]= [0,a_{1},a_{2},\cdots, a_{n}-1, 1]$$ with $a_{n} > 1.$ Then each of the two continued fraction forms for $p/q$ has its own truncate. These truncates are
  $$a/b = [0,a_{1},a_{2},\cdots, a_{n-1}]$$ and 
  $$c/d = [0,a_{1},a_{2},\cdots, a_{n}-1].$$
 The reader will easily prove that for $n$ odd, $a/b < p/q < c/d$ with $aq-bp = pd- cq = ad - bc = -1$
 and that $p/q = (a+c)/(b+d).$ Similarly, if $n$ is even, then $c/d < p/q < a/b,$ with corresponding results.
 Thus the  truncates $a/b$ and $c/d$ form a pair of ancestors for the fraction $p/q$ in the Farey process.
 It then follows by induction that this pair {\it is} the pair from which the fraction $p/q$ appears if the process starts with $0/1$ and $1/1$ as we have outlined it. 
\bigbreak

Here is an example: $4/7 = 1/(1 + 3/4) =  1/(1 + 1/(1 + 1/3)) = [0,1,1,3].$ The truncates are
$[0,1,1] = 1/2$ and $[0,1,1,2] = 3/5.$ We have $$1/2 < 4/7 < 3/5,$$ and indeed $4/7$ arises from
$1/2$ and $3/5$ in the Farey process, as is illustrated above. It is interesting to note that we can then look at the Farey process as a supplier of infinitely many unknot diagrams obtained from these adjacent truncates.
Perhaps the knot theory has deeper connections with number theory than we know.
\bigbreak

\begin{rem} \rm
Recall that the {\it Riemann Hypothesis} says that all the non-trival zeroes of the Riemann zeta function
$$\zeta(s) = \Sigma_{n=1}^{\infty} \frac{1}{n^{s}}$$ lie on the half-line $s = \frac{1}{2} + it$ where $i$ is the square root of $-1$
and $t$ is a real number. This famous unsolved problem is equivalent to a statement about the Farey series:
Let $A(x)$ denote the number of terms in the Farey series corresponding to $x.$
Let $\delta_{j}(x)$ denote the amount by which the $j$-th term of the Farey series for $x$ differs from $j/A(x).$ The following
conjecture is equivalent to the Riemann Hypothesis. 
\begin{conj}{\bf (Franel and Landau).} For each  $\epsilon > 0$ there exists a constant $K$, depending upon $\epsilon,$ such that  
$\Sigma_{j=1}^{A(x)} |\delta_{j}(x)|  < Kx^{\frac{1}{2}+\epsilon}$  as $x \longrightarrow \infty.$
\end{conj}
The equivalence is due to Franel and Landau \cite{FL}. See also \cite{Edwards}.
Thus the Riemann Hypothesis is equivalent to a statement about how certain collections of unknots jostle one another as pairs of rational 
numbers in the real line.
\end{rem}
\bigbreak

In Figure 38 (top) we give an illustration of the Farey generating process for positive rational numbers. Here the tiers of
numbers are connected in a graph, whose top nodes are the ``ancestors" $\frac{0}{1}$ and $\frac{1}{0}.$ A given fraction is
the mediant of its two ancestors. In the graph each fraction has two lines upward to its immediate ancestors. The graph contains a natural
binary tree starting from
$\frac{1}{1}$ (see bottom illustration of Figure 38). We have indicated the genesis of this binary tree as a subgraph of this
graph that shows all ancestors of the Farey fractions. 
\bigbreak

\noindent {\it A given positive continued fraction
$[a_{1},a_{2},...,a_{n}]$ can be found (say $n$ is odd) by starting at $\frac{1}{1}$ in the tree and  heading downward by $a_{1}$ edges to the
right,
$a_{2}$ edges to the left, $a_{3}$ edges to the right, $\ldots$,
$a_{n-1}$ edges to the left, and finally $a_{n} -1$ edges to the right. If $n$ is even, the sequence will start with $a_{1}$ edges to the right
and  end with
$a_{n} -1$ edges to the left.} 
\bigbreak

\noindent For example $\frac{7}{5} = [1,2,2]$ is related to the sequence of
instructions
$RLLR$ where $R$ denotes ``right" and $L$ denotes ``left". The reader will note that these instructions take one from $\frac{1}{1}$ to
$\frac{7}{5}$ in  the tree of Figure 38. For another example, take the instruction $RRLL$ and note that it takes one from
$\frac{1}{1}$ to $\frac{7}{3}.$ We have $\frac{7}{3} = [2,3]$ which indeed is related to $RRLL.$ On the other hand, we also
have $\frac{7}{3} = [2,2,1]$ and this also is related to $RRLL.$ Thus the truncation to $a_{n} - 1$ left or right steps
in the last part of the path corresponds to our earlier discussion of the ambiguity of the last term in the sequence for the continued
fraction.

$$ \picill5inby4.6in(UK38) $$
\begin{center}
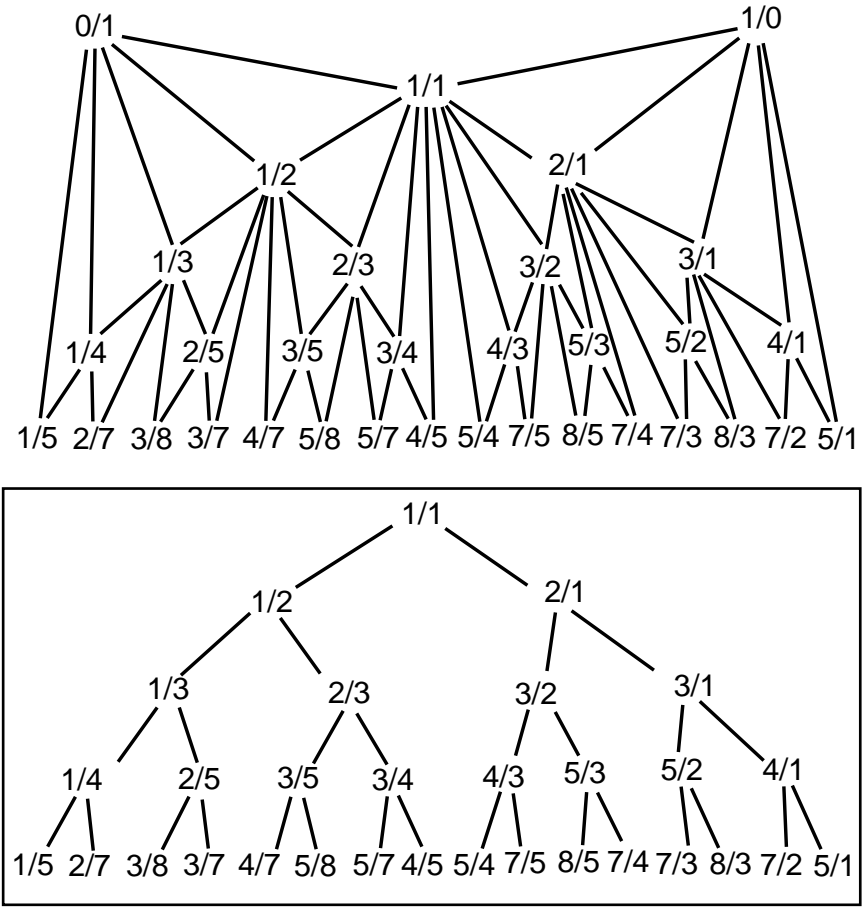
{\bf Figure 38 - The graph of Farey fractions and binary tree of ancestors}
\end{center} 

\noindent The upshot of this structure is that one can use a Farey tree (also called the {\it Stern-Brocot tree}) to enumerate the possible
unknots and their corresponding continued fraction pairs. 
It is easy to see from this prescription that a mediant and one of its ancestors are convergents. Compare these comments with our Theorem 8.
Note also that since the rationals are dense in the real numbers, this tree provides a way to obtain for each positive real number,
a continued fraction that converges to  it. Each real number is the limit of a path going down the binary
tree. For example, the golden ratio $\phi = \frac{1 + \sqrt{5}}{2} = [1,1,1,1,\ldots]$
is the limit of the sequence $\{R,RL,RLR,RLRL,RLRLR,\ldots\}.$ The inverse of the golden ratio, $\frac{1}{\phi} =
\frac{\sqrt{5}-1}{2} = [0,1,1,1,1,\ldots]$, is the limit of the sequence $\{L,LR,LRL,LRLR,LRLRL,\ldots\}.$ Each path down the infinite binary
tree corresponds to the continued fraction expansion of a unique real number.

\subsection {Pick's Theorem}

Pick's Theorem \cite{TheBook}, \cite{Bru} states that the area of a polygon in the standard integral lattice in the Euclidean plane
is given by the formula $$A = I + B/2 -1$$ where $I$ denotes the number of lattice points in the interior of the polygon and $B$ denotes 
the number of lattice points on the boundary. See Figure 39 for an example.

$$ \picill3inby2.3in(UK39) $$
\begin{center}
{\bf Figure 39 - Area of a Lattice Polygon}
\end{center}

\noindent The simplest case of Pick's Theorem is when there are no interior
lattice points and the figure is a triangle with just three lattice points on the boundary. Call such a figure a {\it small triangle.} The
area of a small triangle is, by Pick's Theorem, equal to 
$1/2.$  If one node of the small triangle is at the origin,
then the other two nodes can be viewed as integral vectors
$(a,b)$ and
$(c,d)$ in the plane. One can show that if a small triangle is formed by the origin and the tips of the two vectors, then
the determinant
$ad-bc$ has absolute value equal to $1.$ Hence the two vectors form an alternative basis for the integer lattice in the plane. They satisfy the
formula $|ad - bc| = 1$ and yield unknot fraction pairs $\frac{a}{b}, \frac{c}{d}$ or $\frac{b}{a}, \frac{d}{c}.$ This means that the
fundamental case of Pick's Theorm is directly related to the structure of unknot pairs and to consecutive fractions in the Farey construction. 
\bigbreak

\noindent We can illustrate this relationship with the integer lattice by plotting vectors in a finite lattice. We will say that two vectors are
{\it adjacent} in the integer lattice if together with the origin they form a triangle with no interior lattice points. Starting with a finite
lattice, one can plot all lattice vectors with non-zero slopes less than or equal to one whose coordinates are relatively prime. This gives a
finite collection of vectors and it is easy to see that two such vectors are adjacent, in the sense given above, if and only if one can be
rotated about the origin into the other without encountering another vector in the collection. Thus adjacency becomes rotational adjacency in
such a plot. In Figure 40 we illustrate this pattern by plotting all such vectors in the $5 \times 5$
lattice. This reproduces the Farey
series for denominators no larger than $5$ if we associate the vector $(a,b)$ with the fraction
$b/a.$ The reader will note that the collection of fractions corresponding to the points in the figure are
$$\{1/5,\,\,\,1/4,\,\,\,1/3,\,\,\,2/5,\,\,\,1/2,\,\,\,3/5,\,\,\,3/4,\,\,\,1/1\}.$$
This is the Farey series for $5\frac{1}{2}.$
\bigbreak

The topic of Pick's Theorem brings us back to topology in that the general case of the Theorem follows from 
the special case of the triangle with no interior lattice points by the use of Euler's formula for plane graphs:  $$V-E+F=2$$
when the connected graph in the plane has $V$ nodes, $E$ edges and $F$ faces (including the outer, unbounded, face). One triangulates the
polygon using the lattice points. All the  faces are small triangles except for the outer face. From this it follows that $3(F-1) + B = 2E,$ and
we have that
$V=I+B.$  These two equations and the Euler formula give $(F-1)/2 = I + B/2 - 1.$ The area of the polygon is $(F-1)/2$ since it is composed of
$F-1$ triangles, each of area one-half. This gives Pick's Theorem. It is curious to think of the triangulation of the 
polygon giving us $F-1$ unknots to examine.

$$ \picill3inby3.8in(UK40) $$
\begin{center}
{\bf Figure 40 - Lattice Vectors and Farey Fractions}
\end{center}

\subsection{Ford Circles}

Finally, there is an amazing geometric interpretation of pairs of fractions $a/b$ and $c/d$ such that $|ad-bc|=1.$
Associate to each reduced fraction $a/b$ on the real line a circle tangent to the real line  of diameter $1/b^{2}.$ 
This is the {\it Ford circle} associated with the fraction (See, for example, \cite{Ford},\cite{C2}). One can easily show that {\it two
fractions form an unknot pair if and only if their Ford circles are tangent to one another}. This means that the mediant of an unknot pair
produces a new circle that is tangent to both of its ancestral Ford circles. See Figure 41 for an illustration of this geometry.
\bigbreak

Ford circles can be regarded as curves in the complex plane. The set of Ford circles is invariant under the action of the  modular group of
transformations of the complex plane. By interpreting the upper half of the complex plane as a model of
the hyperbolic plane (the Poincar\'{e} half-plane model) the set of Ford circles can be interpreted as a tiling of the hyperbolic plane. Thus,
any two Ford circles are congruent in hyperbolic geometry. If $C$ and $C'$ are two tangent Ford circles, then the half-circle joining their
respective fractions (on the real line) that is perpendicular to the real line is a hyperbolic line that also passes through the tangent point
of the two circles. This half-circle is indicated with a dashed curve in Figure 41.

$$ \picill3inby2in(UK41) $$
\begin{center}
{\bf Figure 41 - Ford Circles}
\end{center} 

We have seen many aspects of the topology, number theory and geometry related to  unknot pairs in the course of this paper.
We hope that this last section stimulates the curiosity of the reader to make new explorations in this domain.
\bigbreak


 \bigbreak

 \noindent {\sc L.H. Kauffman: Department of Mathematics, Statistics and
 Computer Science, University
 of Illinois at Chicago, 851 South Morgan St., Chicago IL 60607-7045,
 U.S.A.}

 \vspace{.1in}
 \noindent {\sc S. Lambropoulou: Department of Mathematics, 
 National Technical University of Athens,
 Zografou Campus, GR-157 80 Athens, Greece. }

\vspace{.1in}
\noindent {\sc E-mails:} \ {\tt kauffman@math.uic.edu  \ \ \ \ \ \ \ \ sofia@math.ntua.gr 

\noindent http://www.math.uic.edu/$\tilde{~}$kauffman/ \ \ \ \ http://www.math.ntua.gr/$\tilde{~}$sofia}


 \end{document}